\pdfoutput=1
\documentclass[12pt]{amsart}
\usepackage{tikz, stmaryrd}
\usetikzlibrary{matrix,arrows,decorations.pathmorphing}

\usepackage{color}
\usepackage{amsmath,amssymb, amsfonts, amscd}
\usepackage[all]{xy}
\usepackage{verbatim}
\usepackage{enumerate}
\usepackage{yhmath}
\usepackage{xcolor}
\usepackage{hyperref}
\newtheorem{thm}{Theorem}[section]
\newtheorem{cor}[thm]{Corollary}
\newtheorem{lem}[thm]{Lemma}
\newtheorem{obs}[thm]{Observation}

\newtheorem{prop}[thm]{Proposition}

\theoremstyle{definition}

\newtheorem{defn}[thm]{Definition}

\newtheorem{example}[thm]{Example}

\newtheorem{rem}[thm]{Remark}
\numberwithin{equation}{section}

\DeclareFontFamily{U}{rsf}{} \DeclareFontShape{U}{rsf}{m}{n}{
  <5> <6> rsfs5 <7> <8> <9> rsfs7 <10->  rsfs10}{}
\DeclareMathAlphabet{\mathscr}{U}{rsf}{m}{n}

\newcommand{\A}{A}

\newcommand{\D}{D}

\renewcommand{\imath}{\sqrt{-1}}

\newcommand{\Ab}{\mathbb A}

\newcommand{\Nb}{\mathbb N}

\newcommand{\Lb}{\mathbb L}
\newcommand{\Rb}{\mathbb R}
\newcommand{\Cb}{\mathbb C}

\DeclareMathOperator{\Hom}{Hom}

\DeclareMathOperator{\spec}{spec}
\DeclareMathOperator{\id}{e}
\DeclareMathOperator{\coim}{coim}
\DeclareMathOperator{\im}{im}
\DeclareMathOperator{\coker}{coker}

\newcommand{\uHom}{\underline{\Hom}}
\newcommand{\ootimes}{\overline{\otimes}}
\newcommand{\wotimes}{\widehat{\otimes}}

\def\sC{{\mathscr C}}

\def\sL{{\mathscr L}}

\def\sT{{\mathscr T}}

\newcommand{\fm}{\mathfrak{m}}

\newcommand{\tC}{\text{\bfseries\sf{C}}}
\newcommand{\ttsC}{\text{\bfseries\sf{sC}}}
\newcommand{\tD}{\text{\bfseries\sf{D}}}

\newcommand{\tT}{\text{\bfseries\sf{T}}}

\newcommand{\tP}{\text{\bfseries\sf{P}}}

\newcommand{\mM}{\mathcal{M}}

\newcommand{\mH}{\mathcal{H}}

\newcommand{\mV}{\mathcal{V}}
\newcommand{\mW}{\mathcal{W}}

\newcommand{\mT}{\mathcal{T}}
\newcommand{\mO}{\mathcal{O}}
\newcommand{\mB}{\mathcal{B}}
\newcommand{\mA}{\mathcal{A}}
\newcommand{\mD}{\mathcal{D}}

\newcommand{\ttMod}{{\text{\bfseries\sf{Mod}}}}
\newcommand{\ttBorn}{{\text{\bfseries\sf{Born}}}}
\newcommand{\ttSNrm}{{\text{\bfseries\sf{SNrm}}}}
\newcommand{\ttNrm}{{\text{\bfseries\sf{Nrm}}}}
\newcommand{\ttCBorn}{{\text{\bfseries\sf{CBorn}}}}

\newcommand{\ttBan}{{\text{\bfseries\sf{Ban}}}}

\newcommand{\Mod}{{\text{\bfseries\sf{Mod}}}}
\newcommand{\ttVect}{{\text{\bfseries\sf{Vect}}}}
\newcommand{\ttComm}{{\text{\bfseries\sf{Comm}}}}

\newcommand{\ttNBorn}{{\text{\bfseries\sf{NBorn}}}}
\newcommand{\ttNTc}{{\text{\bfseries\sf{NTc}}}}

\newcommand{\ttC}{{\tC}}
\newcommand{\ttD}{{\tD}}
\newcommand{\ttP}{{\tP}}

\newcommand{\ttT}{{\tT}}

\newcommand{\ttA}{{\text{\bfseries\sf{A}}}}

\newcommand{\ttAff}{{\text{\bfseries\sf{Aff}}}}
\newcommand{\ttAfnd}{{\text{\bfseries\sf{Afnd}}}}
\newcommand{\ttStn}{{\text{\bfseries\sf{Stn}}}}

\newcommand{\ttLoc}{{\text{\bfseries\sf{Tc}}}}

\newcommand{\ttInd}{{\text{\bfseries\sf{Ind}}}}
\newcommand{\ttPro}{{\text{\bfseries\sf{Pro}}}}

\newcommand{\ttInt}{{\text{\bfseries\sf{Int}}}}

\newcommand{\ttCpt}{{\text{\bfseries\sf{Cpt}}}}

\newcommand{\ie}{\emph{i.e. }}
\newcommand{\eg}{\emph{e.g. }}
\newcommand{\la}{\lambda}
\renewcommand{\a}{\alpha}

\def\then{\Rightarrow}
\newcommand{\ol}{\overline}

\newcommand{\what}{\widehat}
\def\rhook{{\hookrightarrow}}
\def\lt{\langle}
\def\gt{\rangle} 
\def\diss{{\rm diss\,}}
\def\void{{\rm \varnothing}}
\def\sup{{\rm sup\,}}
\def\Max{{\rm Max\,}}
\def\ho{{\rm ho\,}}

\def\limpro{\mathop{\lim\limits_{\displaystyle\leftarrow}}}
\def\limind{\mathop{\lim\limits_{\displaystyle\rightarrow}}}

\numberwithin{equation}{section}

\begin{document}

\title[]{Stein Domains in Banach Algebraic Geometry}
\author{Federico Bambozzi, Oren Ben-Bassat, and Kobi Kremnizer}\thanks{}%
\address{Federico Bambozzi, Fakult\"{a}t f\"{u}r Mathematik,
Universit\"{a}t Regensburg,
93040 Regensburg, 
Germany}
\bigskip
\email{Federico.Bambozzi@mathematik.uni-regensburg.de}
\address{Oren Ben-Bassat,
Department of Mathematics,
University of Haifa,
Haifa, Israel}
\hfill \break
\bigskip
\email{ben-bassat@math.haifa.ac.il}
\address{Kobi Kremnizer, Mathematical Institute,
Radcliffe Observatory Quarter,
Woodstock Road,
Oxford
OX2 6GG, UK}
\hfill \break
\bigskip
\email{Yakov.Kremnitzer@maths.ox.ac.uk}
\dedicatory{}
\subjclass{}%
\thanks{
        The first author acknowledges the support of the University of Padova by MIUR PRIN2010-11 ``Arithmetic Algebraic Geometry and Number Theory", and the University of Regensburg with the support of the DFG funded CRC 1085 ``Higher Invariants. Interactions between Arithmetic Geometry and Global Analysis" that permitted him to work on this project. The second author acknowledges the University of Oxford and the support of the European Commission under the Marie Curie Programme for the IEF grant STACKSCATS which enabled this research to take place. The contents of this article reflect the views of the authors and not the views of the European Commission. We would like to thank Konstantin Ardakov, Francesco Baldassarri, Elmar Grosse-Kl\"onne, Fr\'{e}d\'{e}ric Paugam, and J\'{e}r\^{o}me Poineau for interesting conversations.}
\keywords{}%

\begin{abstract}
\noindent In this article we give a homological characterization of the topology of Stein spaces over any valued base field. In particular, when working over the field of complex numbers, we obtain a characterization of the usual Euclidean (transcendental) topology of complex analytic spaces. For non-Archimedean base fields the topology we characterize coincides with the topology of the Berkovich analytic space associated to a non-Archimedean Stein algebra. Because the characterization we used is borrowed from a definition in derived geometry, this work should be read as a derived perspective on analytic geometry.
\end{abstract}
\maketitle 
\tableofcontents
\section{Introduction} 

This paper can be thought as a continuation of \cite{BeKr} and \cite{BaBe}, on which we build our main results. The problem addressed in the present work is to extend the approach proposed in \cite{BeKr} and \cite{BaBe} to Stein domains of analytic spaces. More precisely, the main theorems of this paper prove that the characterization of open embeddings with the homotopy monomorphism property, given in \cite{BeKr} for affinoid domains and in \cite{BaBe} for dagger affinoid domains, extends to Stein domains.
Since Stein spaces are defined as spaces which have a suitable exhaustion by (dagger) affinoid subdomains, the strategy of our proofs is to rely on the results of \cite{BeKr} and \cite{BaBe}, and deduce the theorems for Stein spaces using their exhaustions. In particular, this strategy naturally leads to consider projective limits of bornological spaces and the issue about the commutation of the bornological projective tensor product with projective limits. We devote Section \ref{sec:functional_analysis} to provide basic results about these issues for which it seems that there is not much material in literature. Our main difficulty is the fact that neither the projective tensor product nor the injective tensor product of bornological vector spaces commutes with projective limits in general. This is also the main obstacle to the generalization of our results to more general notion of domains, such as dagger quasi-Stein domains.

The paper is organized as follows: in Section \ref{sec:2} we recall the main results from \cite{BaBe} that will be used, we fix the notation and we introduce some key notions that will be used in our main proofs. Section \ref{sec:functional_analysis} contains our technical results in functional analysis that are required to prove Theorems \ref{thm_stein_homotopy} and \ref{thm:DaggerQStHoEpToImm}. In Section \ref{sec:proper} we describe the class of proper bornological spaces and we prove its main property in Proposition \ref{prop:sequence_close}, which shows that the closure of a subspace of a proper bornological space is equal to the set of its (bornological) limit points. Section \ref{sec:nuclear} begins with the definition of nuclear bornological spaces and the study of their main properties. We remark that, although proper and nuclear bornological vector spaces have been previously considered in literature, see for example \cite{Hog} and \cite{Hog3}, this is the first time, at our knowledge, that a detailed study of their property is done also for non-Archimedean base fields. The main result of Section \ref{sec:nuclear} is Theorem \ref{thm:strct_exact_nuclear} which shows that nuclear bornological vector spaces are flat in $\ttCBorn_k$ with respect to its natural monoidal structure, \ie the endofunctor $(-) \wotimes_k F$ is exact if $F$ is nuclear. This is a bornological counterpart of the well-known result of the theory of locally convex spaces. The rest of Section \ref{sec:functional_analysis} deals with projective limits of bornological vector spaces: In particular, Section \ref{sec:relative_flat} addresses the issue of commutation of projective limits with the bornological complete projective tensor product; Section \ref{sec:exact_sequences} deals with a bornological version of the Mittag-Leffler Lemma for Fr\'echet spaces and Section \ref{sec:der_lim} contains some lemmas about the computation of the derived functor of the projective limit functor for quasi-abelian categories. Section \ref{sec:Stein} starts by providing some results on Fr\'echet bornological algebras and then it continues by giving the definition of Stein spaces and Stein algebras suitable for our context. The main result of this section, Theorem \ref{prop:quasi_Stein_algebra_spaces}, is a generalization of Forster's Theorem about the anti-equivalence between the category of complex Stein algebras and complex Stein spaces, to arbitrary base fields. Section \ref{sec:Stein_geometry} contains our main results. We characterize the open embeddings of Stein spaces by the maps having the homotopy monomorphism property (see Theorem \ref{thm_stein_homotopy} and Theorem \ref{thm:DaggerQStHoEpToImm}). Giving a disjoint union of Stein spaces mapping to a fixed Stein space, we characterize in Theorem \ref{thm:coverings} the surjectivity of such a morphism by a conservativity property for transversal modules. See Definition \ref{defn:RRqcoh} for the notion of transversal module, called a RR-quasicoherent module, after the work of Ramis and Ruget \cite{RR}.

The homological treatment of holomorphic functional calculus, as developed by Taylor, leads naturally to derived methods (see for example \cite{Tay}). Let $(\ttC, \ootimes, \id_{\ttC})$ be a closed symmetric monoidal elementary quasi-abelian category with enough flat projectives. In \cite{BeKr2} a Grothendieck topology of (homotopy) Zariski open immersions in $\ttComm(\ttsC)^{op}$ is defined where $\ttsC$ is the closed symmetric monoidal model category of simplicial objects in $\ttC$.  The homotopy monomorphism condition that we use in this article (in the case that $\ttC= \ttCBorn_{k}$ or $\ttInd(\ttBan_{k})$) is a restriction of Definition 1.2.6.1 (3) of \cite{TVe3} from the homotopy category of $\ttComm(\ttsC)^{op}$ to the opposite category of dagger Stein algebras (thought of as constant simplicial objects). The model structure to be explained in \cite{BeKr2} is compatible in a natural way with the quasi-abelian structure on $\ttC$. This allows us to relate our work to the work of To\"{e}n and Vezzosi from \cite{TVe2, TVe3}, as shown in \cite{BeKr2}: $\ttComm(\ttsC)^{op}$ satisfies their axioms on a monoidal model category so according to their approach one can do derived geometry relative to it, and apply their results. 

In the case that our base field is the complex numbers, some of our results are already present in the work of Pirkovskii, for instance \cite{Pir}. 
\subsection{Notation}

The notation used here will be totally consistent with the notation of \cite{BaBe}, which is the following:

\begin{itemize}
	\item If $\ttC$ is a category we will use the notation $X \in \ttC$ to indicate that $X$ is an object of $\ttC$.
	\item If $\ttC$ is a category then $\ttInd(\ttC)$ will denote the category of Ind-objects of $\ttC$.
	\item $k$ will denote a field complete with respect to a fixed non-trivial valuation, Archimedean or non-Archimedean.
	\item $\ttVect_{k}$ is the closed symmetric monoidal category of vector spaces (with no extra structure) over $k$.
	\item $\ttSNrm_k$ the category of semi-normed modules over $k$, remarking that, if not otherwise stated, by a semi-normed space over non-Archimedean base field we mean a $k$-vector space equipped with a non-Archimedean semi-norm.
	\item $\ttNrm_k$ the category of normed modules over $k$.
	\item $\ttBan_k$ the category of Banach modules over $k$.
	\item For $V\in \ttSNrm_k$, $V^{s}= V/\overline{(0)}$ is the separation and $\widehat{V} \in \ttBan_{k}$ is the separated completion.
	\item $\ttBorn_{k}$ the category of bornological vector spaces of convex type over $k$ and $\ttCBorn_{k}$ the category of complete bornological vector spaces of convex type over $k$.
	\item For $E \in \ttBorn_{k}$ and $B$ a bounded absolutely convex subset, (i.e. a bounded disk) of $E$,  $E_B$ is the linear subspace of $E$ spanned by elements of $B$ equipped with the gauge semi-norm (also called the Minkowski functional) defined by $B$ (see Remark 3.40 of \cite{BaBe} for a review of the notion of gauge semi-norm).
	\item  For $E \in \ttBorn_{k}$, $\mathcal{D}_{E}$ denotes the category of bounded absolutely convex subsets of $E$.
	\item $\ttAfnd^{\dagger}_{k}$ denotes the category of dagger affinoid algebras over $k$.
	\item  For $E \in \ttBorn_{k}$, $\mathcal{D}^{c}_{E}$ denotes the category of bounded absolutely convex subsets $B$ of $E$ for which $E_{B}\in \ttBan_{k}$.
	\item The notation $\limind$ refers to a colimit (also known as inductive or direct limit) of some functor in a category.
	\item The notation $\limpro$ refers to a limit (also known as projective or inverse limit) of some functor in a category.
	\item For polyradii $\rho = (\rho_i) \in \Rb_+^n$ and $\rho'$, the notation $\rho < \rho'$ means that $\rho$ and $\rho'$ have the same number of components and every component of $\rho$ is strictly smaller than the corresponding component of $\rho'$. 
	\item With the notation $\ttC^\ttD$ we will denote the category of covariant functors $\ttD \to \ttC$. In particular, if $I$ is a filtered set we will denote $\ttC^I$ the category of functors $I \to \ttC$, when $I$ is thought as a category. 
	\item A cofiltered projective system $\{ E_i \}_{i \in I}$ of objects of a category is said to be \emph{epimorphic} if for any $i < j$ the system map $E_j \to E_i$ is an epimorphism. Similarly, a filtered direct system $\{ E_i \}_{i \in I}$ is said to be \emph{monomorphic} if for any $i < j$ the system map $E_i \to E_j$ is a monomorphism. 
\end{itemize}

\section{Bornological Algebraic Geometry} \label{sec:2}

\subsection{Quasi-abelian categories, bornological spaces and dagger analytic geometry}
We suppose that the reader is familiar with the theory of quasi-abelian categories as developed in \cite{SchneidersQA}. In this section $(\ttC, \ootimes, \id_{\ttC})$ will be a closed symmetric monoidal quasi-abelian category and $\uHom$ will denote the internal hom functor. To any closed symmetric monoidal category is associated a category of commutative monoids, denoted $\ttComm(\ttC)$ and a category of affine schemes $\ttAff(\ttC) = \ttComm(\ttC)^{op}$. The duality functor $\ttComm(\ttC) \to \ttAff(\ttC)$ is denoted by $\spec$.
To any $A \in \ttComm(\ttC)$ we can associate the category of $A$-modules $\ttMod(A)$, which is quasi-abelian closed symmetric monoidal with respect to a bifunctor $\ootimes_A$, naturally induced by $\ootimes$. Moreover, since $\ttMod(A)$ is quasi-abelian we can always associate to $A$ the derived categories of $\ttMod(A)$, denoted $D(A)$, and using the left t-structure of $D(A)$ we define $D^{\leq 0}(A)$, $D^{\geq 0}(A)$ and $D^b(A)$. Notice also that in Proposition 2.1.18 (c) of \cite{SchneidersQA} it is shown that if $\ttC$ is elementary quasi-abelian, then $\ttMod(A)$ is elementary quasi-abelian.

\begin{defn}
	A morphism $\spec(B)\to \spec(A)$ is said to be a \emph{homotopy monomorphism} if the canonical functor $D^{\leq 0}(B)\longrightarrow D^{\leq 0}(A)$ is fully faithful. In a dual way, we say that the correspondent morphism of monoids $A \to B$ is a \emph{homotopy epimorphism}.
\end{defn}

The following characterization of homotopy monomorphisms is the useful one for practical purposes.

\begin{lem} \label{lem_HomotopyMon}
	Assume that $p:\spec(B)\to \spec(A)$ is a morphism in $\ttAff(\ttC)$ and that the functor $\ttMod(A) \to \ttMod(B)$ given by tensoring with $B$ over $A$ is left derivable to a functor $D^{\leq 0}(A)\to D^{\leq 0}(B)$. Then, $p$ is a homotopy monomorphism if and only if $B\ootimes^{\mathbb{L}}_{A} B\cong B$.
\end{lem}

{\bf Proof.} 
See Lemma 2.24 of \cite{BaBe}.
\ \hfill $\Box$

\begin{lem} \label{lem:composition_HomotopyMon}
	Let $f: \spec(A) \to \spec(B)$, $g: \spec(B) \to \spec(C)$ be two morphisms of affine schemes such that $g \circ f$ and $g$ are homotopy monomorphisms, then also $f$ is a homotopy monomorphism.
\end{lem}
{\bf Proof.} 
The hypothesis mean that we have a diagram of functors
\[
\begin{tikzpicture}
\matrix(m)[matrix of math nodes,
row sep=2.6em, column sep=2.8em,
text height=1.5ex, text depth=0.25ex]
{ D^{\leq 0}(A) & & D^{\leq 0}(B)  \\
	& D^{\leq 0}(C) \\};
\path[->,font=\scriptsize]
(m-1-1) edge node[auto] {$f_*$} (m-1-3);
\path[->,font=\scriptsize]
(m-1-1) edge node[auto] {$(g\circ f)_*$} (m-2-2);
\path[<-,font=\scriptsize]
(m-2-2) edge node[auto] {$g_*$}  (m-1-3);
\end{tikzpicture}
\]
such that $g_*$ is fully faithful and $g_* \circ f_*$ is fully faithful. Hence for any $V,W \in D^{\leq 0}(A)$
\[ \Hom_{D^{\leq 0}(A)}(V, W) \cong \Hom_{D^{\leq 0}(C)}((g\circ f)_*(V), (g \circ f)_*(W)) \cong \Hom_{D^{\leq 0}(B)}(f_*(V), f_*(W)) \]
which precisely means that $f_*$ is fully faithful.
\ \hfill $\Box$

We recall the following notion from \cite{BaBe} and \cite{BeKr}.

\begin{defn}\label{defn:RRqcoh}
	Consider an object $A \in \ttComm(\ttCBorn_k)$. We define a sub-category $\ttMod^{RR}(A)$ of $\ttMod(A)$ whose modules $M$ satisfy the property that the natural morphism $M \wotimes^{\mathbb{L}}_{A} B \to M \wotimes_{A}B$ is an isomorphism in $\D^{\le 0} (B)$, for all homotopy epimorphisms $A \to B$. We call these modules \emph{RR-quasicoherent modules}.
\end{defn}

Homotopy epimorphisms are the morphisms that we use to endow $\ttComm(\ttC)$ with a Grothendieck topology. The following definition is based on definitions in \cite{TV} and \cite{TVe3}.

\begin{defn}\label{defn:homotopy_Zariski}
	Consider a full sub-category $\ttA \subset \ttAff(\ttC)$ such that the base change of a homotopy monomorphisms in $\ttA$ is a homotopy monomorphism. On $\ttA$ we can define the \emph{homotopy Zariski topology} which has as its covers collections $\{\spec (B_i)\to \spec(A)\}_{i\in I}$ where there exists a finite subset $J \subset I$ such that
	\begin{itemize}
		\item
		for each $i\in J$, the morphism $A\to B_i$ is of finite presentation and the resulting morphisms $D^{\leq 0}(B_i)\to D^{\leq 0}(A)$ is fully faithful;
		\item
		a morphism in $\ttMod^{RR}(A)$ is an isomorphism if and only if it becomes an isomorphism in each $\ttMod^{RR}(B_j)$ for $j\in J$ after applying the functor $M \mapsto M\ootimes_{A}^{\mathbb{L}} B_{j}$. Such a family is called \emph{conservative}.
	\end{itemize}
	One can drop the requirement on the maps of the covering $A\to B_i$ to be finitely presented, obtaining another topology called the \emph{formal homotopy Zariski topology}. 
\end{defn}

Later on, we will relax the condition on coverings allowing the subset $J \subset I$ to be countable. We will discuss this issue at the end of Section \ref{sec:Stein_geometry}.
We will also make an extensive use of flat objects in $\ttC$, in the following sense.

\begin{defn}\label{defn:flat}Let $(\ttC, \ootimes, \id_{\ttC})$ be a closed, symmetric monoidal, quasi-abelian category. We call an object $F$ of $\ttC$ \emph{flat} if for any strictly exact sequence 
	\[0 \to E' \to E \to E'' \to 0
	\]
	the resulting sequence 
	\[0 \to E'\ootimes F \to E\ootimes F \to E'' \ootimes F \to 0
	\]
	is strictly exact, \ie if the endofunctor $E \mapsto E \ootimes F$ is an exact functor in the terminology of \cite{SchneidersQA}.
\end{defn}

We conclude this section by defining free resolutions in closed symmetric monoidal quasi-abelian categories and by showing some properties of them.

\begin{defn}
	Let $(\ttC, \ootimes, \id_{\ttC})$ be a closed symmetric monoidal quasi-abelian category and let $A \in \ttComm(\ttC)$.  An object $E\in \ttMod(A)$ is called \emph{free} if 
	\[ E \cong A \ootimes V \]
	for some $V \in \ttC$.
\end{defn}
\begin{defn} \label{def:free_resolution}
	Let $(\ttC, \ootimes, \id_{\ttC})$ be a closed symmetric monoidal quasi-abelian category and let $A \in \ttComm(\ttC)$. A \emph{free resolution} of $E \in \ttMod(A)$ is the data of a strict complex 
	\[ \cdots \to L^{2}(E) \to L^{1}(E) \to L^{0}(E) \to 0 \] 
	and a strict quasi-isomorphism 
	\[ L^\bullet(E) \cong E \]
	where each $L^i(E)$ is free in $\ttMod(A)$ and $E$ is thought as a complex concentrated in degree $0$. 
\end{defn}

\begin{lem} \label{lemma:flat_res}
	Let $(\ttC, \ootimes, \id_{\ttC})$ be a closed symmetric monoidal quasi-abelian category with enough projectives. Let $A \in \ttComm(\ttC)$ and $E \in \ttMod(A)$. Then, $E$ admits a free resolution. If in addition, both $E$ and $A$ are flat as objects in $\ttC$ then each term of the free resolution can be chosen to be a flat object in $\ttMod(A)$.
\end{lem}  
{\bf Proof.}
Consider 
\[ \mathscr{L}_{A}^{n}(E) = A \ootimes (\underbrace{A \ootimes \cdots \ootimes A}_{n \text{ times}} \ootimes E) \]
where we think the first $A$ factor as an $A$-module and the other factors as objects of $\ttC$. In this way $\mathscr{L}_{A}^{n}(E)$ is by definition a free  $A$-module.
Defining the differentials $d_n: \mathscr{L}_{A}^{n}(E) \to \mathscr{L}_{A}^{n-1}(E)$ in the following way: Let $m_A: A \ootimes A \to A$ denotes the multiplication map of $A$ and $\rho_E: A \ootimes E \to E$ the action of $A$ on $E$, then
\begin{equation} \label{eqn:diff_bar}
d_n = \sum_{i = 0}^{n-1} (-1)^i id_A \ootimes \cdots \ootimes m_A \ootimes \cdots \ootimes id_A \ootimes id_E + (-1)^n id_A \ootimes \cdots \ootimes id_A \ootimes \rho_E,
\end{equation} 
where $m_A$ is at the $i$th place. Standard computations show that the complex $\mathscr{L}_{A}^\bullet(E)$ is a free resolution of $E$. This complex is a resolution because $\mathscr{L}_{A}^{\bullet}(E)$ has a splitting over $\ttC$ given by the maps
\[\mathscr{L}_{A}^{n}(E) \longleftarrow \mathscr{L}_{A}^{n - 1}(E)
\]
\begin{equation} \label{eqn:splitting}
1_A \ootimes id_{\mathscr{L}_{A}^{n - 1}(E)},
\end{equation}  
where $1_A$ is the constant morphism to the identity of $A$.
Therefore, we can deduce that the cone of the map $\mathscr{L}_{A}^{\bullet}(E) \to E$ is a strictly exact complex in $D^{\leq 0}(\ttC)$. By Proposition 1.5.1 of \cite{SchneidersQA} a morphism in $\Mod(A)$ is strict if and only if is strict as a morphism in $\ttC$, hence the cone $\mathscr{L}_{A}^{\bullet}(E)  \to E$ is also a strictly exact complex in $D^{\leq 0}(A)$.

It remains to show the claim about the flatness of $\sL_A^n(E)$. Since $\ootimes_A$ is right exact we only need to show only the left exactness of the functor $(-) \ootimes_A \sL_A^n(E)$. Consider a strictly exact sequence of morphisms of $\ttMod(A)$
\[ 0 \to F \to G \to H. \]
Applying $(-) \ootimes_A \sL_A^n(E)$ we obtain the sequence
\begin{equation} \label{eqn:flat_bar_res}
0 \to F \ootimes_A \sL_A^n(E) \to G \ootimes_A \sL_A^n(E) \to H \ootimes_A \sL_A^n(E)
\end{equation}
which can be rewritten as
\[ 0 \to F \ootimes (\underbrace{A \ootimes \cdots \ootimes A}_{n \text{ times}} \ootimes E) \to G \ootimes (\underbrace{A \ootimes \cdots \ootimes A}_{n \text{ times}} \ootimes E) \to H \ootimes (\underbrace{A \ootimes \cdots \ootimes A}_{n \text{ times}} \ootimes E). \]
Therefore, the hypothesis that $A$ and $E$ are flat objects of $\ttC$ directly implies that the sequence (\ref{eqn:flat_bar_res}) is strictly exact.
\ \hfill $\Box$

\begin{defn}\label{defn:Bar}
Let $(\ttC, \ootimes, \id_{\ttC})$ be a closed symmetric monoidal quasi-abelian category. Let $A \in \ttComm(\ttC)$ and $E \in \ttMod(A)$. We define the \emph{Bar resolution} of $E$ to be the free resolution introduced in Lemma \ref{lemma:flat_res}.
\end{defn}

Using the Bar resolution we can prove the following important lemma.

\begin{lem} \label{lem:ind_limit_homotopy_epi}
	Let $\ttC$ be elementary quasi-abelian. Let $\{A_i\}_{i \in I}$ be a filtered inductive system in $\ttComm(\ttC)$ such that all system morphisms are homotopy epimorphisms. Then, for any $j \in I$ the canonical maps
	\[ A_j \to \limind_{i \in I} A_i \]
	are homotopy epimorphisms.
\end{lem}

{\bf Proof.} 
Let's fix a $j \in I$. To show that $A_j \to \underset{i \in I}\limind A_i$ is a homotopy epimorphism we need to check that
\[ (\limind_{i \in I} A_i) \ootimes_{A_j}^\Lb (\limind_{i \in I} A_i )\cong \limind_{i \in I} A_i \]
in $D^{\le 0}(\ttC)$. Consider the Bar resolution $\mathscr{L}_{A_j}^{\bullet}(\limind_{i \in I} A_i)$. Then, the complex
\[ (\limind_{i \in I} A_i) \ootimes_{A_j} \mathscr{L}_{A_j}^{\bullet}(\underset{i \in I}\limind A_i) \]
is a representative of  $(\underset{i \in I}\limind A_i )\ootimes_{A_j}^\Lb( \underset{i \in I}\limind A_i)$.
More explicitly, for each $n \in \Nb$
\[ (\limind_{i \in I} A_i) \ootimes_{A_j} \mathscr{L}_{A_j}^{n}(\underset{i \in I}\limind A_i) = (\limind_{i \in I} A_i) \ootimes_{A_j} A_j \ootimes (\underbrace{A_j \ootimes \cdots \ootimes A_j}_{n \text{ times}} \ootimes (\limind_{i \in I} A_i)) \]
which simplifies to
\[ (\limind_{i \in I} A_i) \ootimes (\underbrace{A_j \ootimes \cdots \ootimes A_j}_{n \text{ times}} \ootimes (\limind_{i \in I} A_i)) \cong \limind_{i \in I} (A_i \ootimes \underbrace{A_j \ootimes \cdots \ootimes A_j}_{n \text{ times}} \ootimes A_i). \]
Now, because by Proposition 2.1.16 (c) of \cite{SchneidersQA} $\underset{i \in I}\limind$ is an exact functor when $I$ is filtered we have $\underset{i \in I}\limind \cong \mathbb{L}\underset{i \in I}\limind$. Therefore, 
\[( \limind_{i \in I} A_i )\ootimes_{A_j}^\Lb (\limind_{i \in I} A_i) \cong (\limind_{i \in I} A_i) \ootimes_{A_j} \mathscr{L}_{A_j}^{\bullet}(\underset{i \in I}\limind A_i) \cong \limind_{i \in I}(A_i \ootimes_{A_j} \mathscr{L}_{A_j}^{\bullet}(A_i)) \cong \limind_{i \in I}(A_i \ootimes_{A_j}^\Lb A_i) \cong \limind_{i \in I} A_i. \]
\ \hfill $\Box$

We also introduce the following complex needed for computing \v{C}ech cohomology in the theories we will develop.

\begin{defn}
Given a collection of morphisms $\mathfrak{U}=\{A\to B_{i}\}_{i \in I}$ in $\ttComm(\ttC)$ and $M \in \ttMod(A)$ we have the \v{C}ech-Amitsur complex $\sC_{A}^\bullet(M,\mathfrak{U})$
\[ \prod_{i} (M \ootimes_{A} B_{i}) \to \prod_{i,j} (M \ootimes_{A} B_{i}\ootimes_{A} B_{j}) \to \cdots
\]
and its augmented version which we call the Tate complex $\sT_{A}^\bullet(M,\mathfrak{U})$
\[0 \to M \to \prod_{i} (M \ootimes_{A} B_{i}) \to \prod_{i,j} (M \ootimes_{A} B_{i}\ootimes_{A} B_{j}) \to \cdots
\]
where we use the degree convention 
\[\mathscr{C}^{d}_{A}(M,\mathfrak{U})= \mathscr{T}^{d}_{A}(M,\mathfrak{U})= \prod_{i_1,\dots, i_d} (M \ootimes_{A} B_{i_1}\ootimes_{A} \cdots \ootimes_{A} B_{i_d})\] 
for $d \geq 1$.
\end{defn}

\subsection{Bornological vector spaces}
In this section we recall very briefly the theory of bornological vector spaces. We refer the reader to Section 3.3 of \cite{BaBe} for more details.

\begin{defn}
	Let $X$ be a set. A \emph{bornology} on $X$ is a collection $\mB$ of subsets of $X$ such that
	\begin{enumerate}
		\item $\mB$ is a covering of $X$, \ie $\forall x \in X, \exists B \in \mB$ such that $x \in \mB$;
		\item $\mB$ is stable under inclusions, \ie $A \subset B \in \mB \then A \in \mB$;
		\item $\mB$ is stable under finite unions, \ie for each $n \in \Nb$ and $B_1, ..., B_n \in \mB$, $\underset{i = 1}{\overset{n}\bigcup} B_i \in \mB$.
	\end{enumerate}
	
	The pair $(X, \mB)$ is called a \emph{bornological set}, and the elements of $\mB$ are called \emph{bounded subsets} of $X$ 
	(with respect to $\mB$, if it is needed to specify).  
	A family of subsets $\mA \subset \mB$ is called a \emph{basis} for $\mB$ if for any $B \in \mB$ there exist $A_1, \dots, A_n \in \mA$ such that $B \subset A_1 \cup \dots \cup A_n$.
	A \emph{morphism} of bornological sets $\varphi: (X, \mB_X) \to (Y, \mB_Y)$ is defined to be a bounded map $\varphi: X \to Y$, \ie a map of sets such that 
	$\varphi(B) \in \mB_Y$ for all $B \in \mB_X$. 
\end{defn}

\begin{defn}
	A \emph{bornological vector space} over $k$ is a $k$-vector space $E$ along with a bornology on the underlying set of $E$ for which the maps $(\la, x) \mapsto \la x$ and $(x, y) \mapsto x + y$
	are bounded.
\end{defn}

\begin{defn}
	A bornological vector space is said to be \emph{of convex type} if it has a basis made of absolutely convex subsets. We will denote by $\ttBorn_k$ the category whose objects are the  bornological vector spaces of convex type and whose morphisms are bounded linear maps between them.
\end{defn}

\begin{rem} \label{rem:gauge}
	For every bornological vector space of convex type $E$ there is an isomorphism
	\[ E \cong \limind_{B \in \mD_E} E_B  \]
	where $B$ varies over the family of bounded absolutely convex subsets of $E$ and $E_B$ is the vector subspace of $E$ spanned by elements of $B$ equipped with the gauge semi-norm (also called Minkowski functional) defined by $B$ and each $E_B$ is equipped with the bornology induced by this semi-norm.
\end{rem}

Reasoning as last remark, one can show that there is a functor $\diss: \ttBorn_k \to \ttInd(\ttSNrm_k)$ which is fully faithful, which commutes with all projective limits and direct sums and whose essential image is the sub-category of essential monomorphic objects of $\ttInd(\ttSNrm_k)$, \ie by ind-objects isomorphic to systems which have monomorphisms as system maps. The functor $\diss$ does not commute with cokernels in general.

\begin{defn} \label{defn:separated_born}
	A bornological vector space over $k$ is said to be \emph{separated} if its only bounded linear subspace is the trivial subspace $\{0\}$.
\end{defn}

\begin{defn} \label{defn:complete_born}
	A bornological vector space $E$ over $k$ is said to be \emph{complete} if there is a small filtered category $I$, a functor $I \to \ttBan_{k}$ and an isomorphism
	\[ E \cong \limind_{i \in I} E_i \]
	for a filtered colimit of Banach spaces over $k$ for which the system morphisms are all injective and the colimit is calculated in the category $\ttBorn_k$.
\end{defn}

Let $\ttCBorn_k$ be the full subcategory of $\ttBorn_k$ consisting of complete bornological spaces. Collecting together several facts from the theory of bornological vector spaces. See Section 3.3 of \cite{BaBe} for a more detailed account.

\begin{lem}\label{lem:CBornProps}
 $\ttCBorn_k$ has all limits and colimits and the inclusion functor $\ttCBorn_k \to \ttBorn_k$ commutes with projective limits. Both $\ttBorn_k$ and $\ttCBorn_k$ are closed symmetric monoidal elementary quasi-abelian categories with enough projectives.
\end{lem}

The monoidal structure of $\ttBorn_k$ is denoted by $\otimes_{\pi, k}$ and is called the \emph{projective tensor product}. The one of $\ttCBorn_k$ is $\wotimes_{\pi, k}$ (often written just as $\wotimes_{k}$) and it is called the \emph{completed projective tensor product}. We conclude by noticing the fact that $\ttCBorn_k$ has enough projectives implies that we can always derive right exact functors from $\ttCBorn_k$ to a quasi-abelian category, as happens in the theory of closed symmetric abelian categories when they have enough projective objects.

\subsection{Tensor product of unbounded complexes} 

In this section we extend to the quasi-abelian settings some results of \cite{BN} on the derived functor of the tensor product for abelian categories. We see how the derived functor $(-)\ootimes^\Lb(-): D^{\le 0}(\ttC) \times D^{\le 0}(\ttC) \to D^{\le 0}(\ttC)$ (discussed extensively in \cite{BaBe}) extends to a functor $(-)\ootimes^\Lb(-): D(\ttC) \times D(\ttC) \to D(\ttC)$. Since we are not looking for the utmost generality of the result we suppose in this section that $\ttC$ is an elementary quasi-abelian closed symmetric monoidal category, although it would be sufficient to suppose only that $\ttC$ has exact direct sums and enough projectives.

Recall that the hypothesis of $\ttC$ being elementary implies that $\ttC$ has enough projectives, cf. Proposition 2.1.15 (c) of \cite{SchneidersQA}.

\begin{lem} \label{lem:homotopy_colim}
Let $X \in D(\ttC)$, then the morphism
\[ \ho\limind_{n \in \Nb} \tau^{\le n}(X) \to \limind_{n \in \Nb} \tau^{\le n}(X) \]
is an isomorphism.
\end{lem}
{\bf Proof.}
We can apply the dual statement of Remark 2.3 of \cite{BN} because by Proposition 2.1.15 (a) of \cite{SchneidersQA}, $LH(\ttC)$ satisfies the axiom AB4${}^*$.
\ \hfill $\Box$

\begin{lem} \label{lem:homotopy_iso}
	Every object in $D(\ttC)$ is strictly quasi-isomorphic to a complex of projective objects.
\end{lem}

{\bf Proof.}
Let $X \in D(\ttC)$. Consider the truncated complex $\tau^{\le n}(X)$ (recall that we are using the left t-structure of $D(\ttC)$). There is a canonical map $\tau^{\le n}(X) \to X$. Since $\ttC$ has enough projectives we can always find strict quasi-isomorphisms $P^{\le n} \to \tau^{\le n}(X)$, where $P^{\le n}$ is a complex of projectives which is zero in degree $> n$. The diagram
\[
\begin{tikzpicture}
\matrix(m)[matrix of math nodes,
row sep=2.6em, column sep=2.8em,
text height=1.5ex, text depth=0.25ex]
{   P^{\le n} & P^{\le n + 1}   \\
	\tau^{\le n}(X) & \tau^{\le n + 1}(X)  \\};
\path[->,font=\scriptsize]
(m-2-1) edge node[auto] {$$} (m-2-2);
\path[->,font=\scriptsize]
(m-1-1) edge node[auto] {$$} (m-2-1);
\path[->,font=\scriptsize]
(m-1-2) edge node[auto] {$$}  (m-2-2);
\end{tikzpicture}
\]
defines a morphism $P^{\le n} \to P^{\le n + 1}$ in $D(\ttC)$, because $P^{\le n + 1} \to \tau^{\le n + 1}(X^\bullet)$ is a strict quasi-isomorphism and this morphism can be realized as a morphism of complexes because $P^{\le n}$ and $P^{\le n + 1}$ are made of projectives. So, there is a sequence of morphisms
\[ \ho\limind_{n \in \Nb} \tau^{\le n}(P^{\le n}) \to \ho\limind_{n \in \Nb} \tau^{\le n}(X) \to \limind_{n \in \Nb} \tau^{\le n}(X)  \]
where $\underset{n \in \Nb}\limind \tau^{\le n}(X) \cong X$ and the first map is an isomorphism because is a homotopy colimit of strict quasi-isomorphisms. The second map is an isomorphism by Lemma \ref{lem:homotopy_colim}. Therefore, $X \cong \ho\underset{n \in \Nb}\limind \tau^{\le n}(P^{\le n})$ which by construction is a complex of projectives. 
\ \hfill $\Box$

\begin{defn}
Let $\ttT$ be a triangulated category with direct sums. We say that subcategory of $\ttT$ is \emph{localizing} if it is closed under direct sums.
\end{defn}

We use the following notation: $K(\ttC)$ the homotopy category of $\ttC$ and $K(\ttP)$ the smallest localizing subcategory of $K(\ttC)$ containing the complexes made of projectives.

\begin{rem}
Notice that since $\ttC$ is elementary quasi-abelian, then both $K(\ttC)$ and $D(\ttC)$ have all the direct sums. Indeed, the claim on $K(\ttC)$ is trivial and for the one on $D(\ttC)$, using Proposition 2.1.12 of \cite{SchneidersQA} one deduces that $\ttC$ is derived equivalent to an elementary abelian category and therefore $D(\ttC)$ has all direct sums as consequence of Corollary 1.7 of \cite{BN}.
\end{rem}

\begin{lem} \label{lem:homotopy_cat}
	The composition functor
	\[ K(\ttP) \rhook K(\ttC) \to D(\ttC)  \]
	is an equivalence.
\end{lem}
{\bf Proof.}
Since $K(\ttP)$ is thick subcategory of $K(\ttC)$, we can use the dual statements of Lemmata 2.8-11 of \cite{BN} to deduce that $K(\ttP)$ is a localizing subcategory of $K(\ttC)$. Using Lemma \ref{lem:homotopy_iso} we see that all objects of $D(\ttC)$ are isomorphic to objects of $K(\ttP)$.
\ \hfill $\Box$

\begin{thm} \label{thm:extend_tensor}
	The tensor product can be derived to a functor
	\[ (-)\ootimes^\Lb(-): D(\ttC) \times D(\ttC) \to D(\ttC). \] The restriction to $D^{\le 0}(\ttC) \times D^{\le 0}(\ttC)$ agrees with the derived tensor product based on projective resolutions.
\end{thm}

{\bf Proof.}
The tensor product $(-)\ootimes(-)$ extends to a functor $K(\ttC) \times K(\ttC) \to K(\ttC)$. By Lemma \ref{lem:homotopy_cat}, under our hypothesis, $D(\ttC)$ is equivalent to a subcategory of $K(\ttC)$. Therefore, we can define $(-)\ootimes^\Lb(-)$ as the restriction to $K(\ttP)$ of the extension of $(-)\ootimes(-)$ to $K(\ttC)$.
\ \hfill $\Box$

\section{Some results in functional analysis} \label{sec:functional_analysis}

This section is devoted to prove results in the theory of bornological vector spaces needed in the rest of the paper.

\subsection{Proper bornological vector spaces}  \label{sec:proper}

\begin{defn}
	Let $E$ be a bornological $k$-vector space and $\{x_n\}$ a sequence of elements of $E$. We say that $\{x_n\}_{n \in \Nb}$ \emph{converges (bornologically) to $0$} if there exists a bounded subset $B \subset E$ such that for every $\la \in k^\times$ there exists an $n = n(\la)$ for which
	\[ x_m \in \la B, \forall m > n. \] 
	We say that $\{ x_n \}_{n \in \Nb}$ converges  (bornologically) to $a \in E$ if $\{x_n - a\}_{n \in \Nb}$ converges (bornologically) to zero.
\end{defn}
The idea of this definition is due to Mackey and it is sometimes called Mackey convergence.

\begin{defn}
	Let $E$ be a bornological vector space over $k$.
	\begin{itemize}
		\item a sequence $\{x_n\}_{n \in \Nb} \subset E$ is called \emph{Cauchy-Mackey} if the double sequence $\{x_n - x_m\}_{n,m \in \Nb}$ converges to zero;
		\item a subset $U \subset E$ is called \emph{(bornologically) closed} if every sequence of elements of $U$ which is bornologically convergent in $X$ converges bornologically to an element of $U$. 
	\end{itemize}
\end{defn}

\begin{defn}
	A bornological vector space is called \emph{semi-complete} if every Cauchy-Mackey sequence is convergent.
\end{defn}

The notion of semi-completeness is not as useful as the notion of completeness in the theory of topological vector spaces. We remark that any complete bornological vector space is semi-complete, but the converse is false.

\begin{rem} \label{rem:born_conv}
	The notion of bornological convergence on a bornological vector space of convex type $E = \underset{B \in \mD_E}\limind E_B$, where $\mD_E$ denotes the family of bounded disks of $E$, can be restated in the following way: $\{ x_n \}_{n \in \Nb}$ is convergent to zero in the sense of Mackey if and only if there exists a $B \in \mD_E$ and $N \in \Nb$ such that for all $n > N$, $x_n \in E_B$ and $x_n \to 0$ in $E_B$ (equipped with the semi-norm induced by $B$).
\end{rem}

\begin{rem}\label{rem:BornTop}
The notion of bornologically closed subset induces a topology on $E$, but this topology is neither a linear topology nor group topology in general. However, an arbitrary intersection of bornological closed subsets of a bornological vector space is bornologically closed. 
\end{rem}

From Remark \ref{rem:BornTop} follows that the following definition is well posed.

\begin{defn}
	Let $U \subset E$ be a subset of a bornological vector space. The (bornological) closure of $U$ is the smallest bornologically closed subset of $E$ in which $U$ is contained. We denote the closure of $U$ by $\ol{U}$.
\end{defn}

\begin{defn} \label{defn:born_dense}
	Let $E$ be a bornological vector space over $k$. We will say that a subset $U \subset E$ is \emph{bornologically dense} if the bornological closure of $U$ is equal to $E$. 
\end{defn}

\begin{prop} \label{prop:strict_morphisms_born_close}
	Let $f: E \to F$ be a morphism in $\ttCBorn_k$, then
	\begin{itemize}
		\item $f$ is a monomorphism if and only if it is injective
		\item $f$ is an epimorphism if and only if $f(E)$ is bornologically dense in $F$
		\item $f$ is a strict epimorphism if and only if it is surjective and $F$ is endowed with the quotient bornology; 
		\item $f$ is a strict monomorphism if and only if it is injective, the bornology on $E$ agrees with the induced bornology from $F$ and $f(E)$ is a bornologically closed subspace of $F$.
	\end{itemize}
\end{prop}

{\bf Proof.}
See Proposition 5.6 (a) of \cite{PrSc}.
\ \hfill $\Box$

\begin{defn} \label{def:proper_born}
	A bornological vector space is called \emph{proper} if its bornology has a basis of bornologically closed subsets.
\end{defn}

\begin{rem} \label{rem:proper}
\begin{enumerate}
\item All bornological vector spaces considered in \cite{Bam} and \cite{BaBe} are proper. 
\item Let $E$ be a separated bornological vector space (\ie $\{0\}$ is a closed subset of $E$), then the morphism $E \to \widehat{E}$ is injective (cf. Proposition 17 page 113 of \cite{H2}). Indeed, in this case, in $E \cong \underset{B \in \mD_E}\limind E_B$ then
\[ \widehat{E} = \widehat{\underset{B \in \mD_E}\limind E_B} \cong \underset{B \in \mD_E}\limind \widehat{E_B}. \]
\end{enumerate}
\end{rem}

The main drawback, for our scope, of proper bornological vector spaces is that although they form a closed symmetric monoidal category, this category is not quasi-abelian. We will also need the following property of proper bornological vector spaces.

\begin{prop} \label{prop:proper_proj}
	Any projective limit of proper objects in $\ttBorn_{k}$ is proper.
\end{prop}

{\bf Proof.}
\cite{H2}, Proposition 12, page 113.
\ \hfill $\Box$

\begin{lem} \label{cor:proper_subspace}
	Any subspace of proper bornological vector space, endowed with the induced bornology, is proper.
\end{lem}

{\bf Proof.}
Let $E$ be a proper bornological vector space and let $\{B_i\}_{i \in I}$ be a basis for its bornology made of bornologically closed bounded subsets. Given a subspace $F \subset E$, the family $\{B_i \cap F\}_{i \in I}$ is a basis for the bornology that $E$ induces on $F$. Consider a sequence of elements $\{x_n\}_{n \in \Nb} \subset B_i \cap F$ converging bornologically to $x \in F$ then, since the inclusion map $F \to E$ is bounded, it preserves bornological limits. Therefore, $\{x_n\}_{n \in \Nb}$ converges to $x$ also as a sequence of elements of $E$ and since $B_i$ is closed $x \in B_i$. It follows that $x \in F \cap B_i$ and therefore the claim.
\ \hfill $\Box$

We have to warn the reader that the closure of a subset $X \subset E$ of a bornological vector space is rarely equal to the limit points of convergent sequences of elements of $X$. So, we introduce the following notation.

\begin{defn} \label{def:limit_points}
	Let $E$ be a bornological vector space and $X \subset E$, then $X^{(1)}$ denotes the elements of $E$ that are bornological limits for sequences of elements of $X$.
\end{defn}

Thus, we have the (in general strict) inclusion $X^{(1)} \subset \ol{X}$. We now show that this inclusion is an equality for an important special case: linear subspaces of complete, proper bornological spaces. We need a preliminary lemma.

\begin{lem}  \label{lemma:closed_complete}
	Let $F \subset E$ be a linear subspace of a complete bornological vector space over $k$. Then, $F$ is (bornologically) closed if and only if $F$ is complete.
\end{lem}
{\bf Proof.}
This is a classical result in the theory of bornological vector spaces. For the convenience of the reader we reproduce here the proof that can be found in \cite{Hog}, chapter 3, Proposition 3.2.1. Moreover our proof deals also with the non-Archimedean case of the lemma not covered in \cite{Hog}. 

Suppose that $F \subset E$ is a complete subspace. Let $\{x_n\}_{n \in \Nb}$ be a sequence of elements of $F$ that converges bornologically to $x \in E$. By Remark \ref{rem:born_conv} there exists a bounded disk $B \subset E$ such that $x_n \to x$ in $E_B$. Since $B \cap F$ is a bounded disk in $F$ and $F$ is complete can find a disk $B' \subset F$ such that $B \subset B'$ and $F_{B'}$ is Banach. This shows that $\{x_n\}_{n \in \Nb}$ is a Cauchy sequence in $F_{B'}$, because the map $F_B \to F_{B'}$ is by construction bounded. Hence, the limit of $\{x_n\}_{n \in \Nb}$ is an elements of $F$, which is hence closed.

For the converse, suppose that $F \subset E$ is a closed subspace. Let $\mD^{c}_{E}$ denotes the family of completant bounded disks of $E$, \ie bounded disks $B$ for which $E_B$ is Banach. It is enough to show that for any $B \in \mD^{c}_{E}$ the bounded disk $B \cap F$ of $F$ is completant. So, let $\{ x_n \}_{n \in \Nb}$ be a Cauchy sequence in $F_{B \cap F}$. This implies that $\{ x_n \}_{n \in \Nb}$ is a Cauchy sequence in $E_B$ and therefore it must converge to a point of $E_B$. But $F_{B \cap F}$ is a closed subspace of $E_B$, since it is the intersection of $E_B \subset E$ with a closed subspace of $E$, hence the limit of $\{ x_n \}_{n \in \Nb}$ is in $F_{B \cap F}$. This proves that $F_{B \cap F}$ is a $k$-Banach space and concludes the proof.
\ \hfill $\Box$

\begin{prop} \label{prop:sequence_close}
	Let $F \subset E$ be a subspace of a complete, proper bornological vector space over $k$, then $F^{(1)} = \ol{F}$, where $F^{(1)}$ is as in Definition \ref{def:limit_points}.
\end{prop}
{\bf Proof.}
Let $E = \underset{i \in I}\limind E_i$ be a presentation of $E$ as a direct limit of Banach spaces over $k$ as in Definition \ref{defn:complete_born}. The bornology on $F$ induced by the inclusion $F \subset E$ can be described by 
\[ F \cong \limind_{i \in I} (F \cap E_i), \] 
where each $F \cap E_i$ is endowed by the norm induced by $E_i$.
By Lemma \ref{cor:proper_subspace} $F$ endowed with this bornology is a proper bornological vector space which is easily seen to be separated. By Remark \ref{rem:proper} (2), the inductive system $\{\widehat{F \cap E_i} \}_{i \in I}$ obtained by applying the completion functor of Banach spaces, is monomorphic and hence it defines a complete and proper bornological vector space $\widehat{F} = \underset{i \in I}\limind \widehat{F \cap E_i}$ in which $F$ embeds. Since $\widehat{F}$ is the completion of $F$ and the inclusion $F \to E$ must factor through a map
\[ \phi: \widehat{F} \to E. \]
Moroever, $\phi$ is a strict monomorphism, because for each $i \in I$ the maps
\[ F \cap E_i \to \widehat{F \cap E_i} \to E_i \]
are strict, since all spaces are endowed with the restriction of the same norm, by construction. $\widehat{F}$ is a closed linear subspace of $E$ because it is complete and hence we can apply Lemma \ref{lemma:closed_complete}. This implies that $\ol{F} \subset \im(\phi) \cong \what{F}$. Finally, let $x \in \im(\phi)$ then $x = \phi(y)$ for some $y \in \widehat{F \cap E_i}$ and for some $i \in I$. So, there exists a sequence of elements $\{x_n\}_{n \in \Nb} \subset F \cap E_i$ which converges in norm to $x$, which shows that $ \what{F} \cong \im(\phi) \subset F^{(1)} \subset \ol{F}$.
\ \hfill $\Box$

We remark that in Proposition \ref{prop:sequence_close} the hypothesis that $F \subset E$ is a (linear) subspace is crucial. It is possible to construct subsets $S\subset E$ for which $S^{(1)} \subset \ol{S}$ strictly even when $E$ is a proper object of $\ttCBorn_{k}$ (see \cite{Hog2} Theorem 1, Section II.7).

\begin{rem}
Proposition \ref{prop:sequence_close} in the case $k =\Cb$ was proven in \cite{M2}, cf. Proposition 4.14. The arguments of \cite{M2} are similar to ours, but the terminology and the context are very different.
\end{rem}

\subsection{Relations between bornological and topological vector spaces} \label{sec:normal}

We recall the definitions of two functors from \cite{H2}, $(-)^{t}: \ttBorn_k \to \ttLoc_k$ and $(-)^b: \ttLoc_k \to \ttBorn_k$, where $\ttLoc_k$ is the category of locally convex topological vector spaces over $k$. To a bornological vector space $E$ we associate the topological vector space $E^t$ in the following way: we equip the underlying vector space of $E$ with a topology for which a basis of $0$-neighborhoods is given by \emph{bornivorous subsets}, \ie  subsets that absorb all bounded subsets of $E$.  If $E$ is a locally convex space, $E^b$ is defined to be the bornological vector space obtained by equipping the underlying vector space of $E$ with the \emph{von Neumann} (also called canonical) bornology, whose bounded subsets are the subsets of $E$ absorbed by all $0$-neighborhoods. 
In chapter 1 of \cite{H2} one can find details of these constructions and their properties of which the main one is that \begin{equation}\label{eqn:tbadj} 
(-)^t : \ttBorn_k \leftrightarrows \ttLoc_k : (-)^b
\end{equation}
is an adjoint pair of functors. 

\begin{lem}\label{lem:GivesEquiv}
	It is shown in \cite{H2} (cf. Chapter 2, \S 2 $n^\circ$ 6) that there are full sub-categories of $\ttBorn_k$  and $\ttLoc_k$ for which these functors give an equivalence of categories. 
\end{lem}
\begin{defn}\label{defn:normal_space}
	We call the objects of the categories mentioned in Lemma \ref{lem:GivesEquiv} normal bornological vector spaces or normal locally convex spaces (depending on the ambient category we are thinking them in). We call the corresponding full sub-categories of normal objects $\ttNBorn_{k} \subset \ttBorn_{k}$ and $\ttNTc_{k} \subset \ttLoc_{k}$.
\end{defn}
For two elements $E, F \in \ttNTc_k$ (or $E, F \in \ttNBorn_k$) the notion of boundedness and continuity of a linear map $f: E \to F$ is equivalent. 
\begin{rem} \label{rem:closed_top_bor_nor}
Any normal bornological vector space is proper (cf. Definition \ref{def:proper_born}). One can show this fact by first noticing that in any topological vector space the closure of bounded subset is bounded. Hence, the von Neumann bornology of a locally convex topological vector space always has a basis made of closed (for the topology) subsets. Moreover, every subset which is closed for the topology is also bornologically closed for the von Neumann bornology. For more details on this topic see \cite{H2}, Chapter 2, \S 3, $n^\circ$ 6 and 7.
\end{rem}

\begin{obs}\label{obs:tbequiv}
	The functor ${(-)}^b: \ttLoc_k \to \ttBorn_k$ commutes with projective limits and ${(-)}^t: \ttBorn_k \to \ttLoc_k$ commutes with inductive limits. 
\end{obs}

\begin{obs}\label{obs:Dense2Dense}Let $E$ be a bornological vector space over $k$. It follows that a subset $U \subset E$ is bornologically dense then $U \subset E^{t}$ is topologically dense. Indeed, all bornologically closed subsets for the von Neumann bornology of $E^{t}$ must be bornologically closed for the bornology of $E$ and, as said in Remark \ref{rem:closed_top_bor_nor}, all topologically closed subsets of $E^{t}$ are bornologically closed for the von Neumann bornology. Therefore, the closure of $U$ in $E^{t}$ contains the closure of $U$ in $E$.
\end{obs}

\begin{example}\label{example:normal}
	\begin{itemize}
		\item All semi-normed spaces over $k$ are in an obvious way normal objects in $\ttBorn_{k}$ and $\ttLoc_{k}$.
		\item All metrizable locally convex vector spaces and metrizable bornological vector spaces of convex type (when equipped with the von Neumann bornology) are normal. See the beginning of page 109 of \cite{H2} or Proposition 3 at page 50 of \cite{Hog}.
		\item The underlying object in $\ttCBorn_{k}$ of a $k$-dagger affinoid algebra is a normal bornological vector space (see \cite{Bam}, Chapter 3).
	\end{itemize}
\end{example}

We introduce some classes of normal spaces that will be used throughout the paper.

\begin{defn} \label{defn:frechet}
A \emph{Fr\'echet space} over $k$  is a complete, metric, locally convex topological vector space over $k$. A \emph{bornological Fr\'echet space} over $k$ is a bornological vector space $E \in \ttCBorn_k$ such that $E \cong F^b$ for a Fr\'echet space $F \in \ttLoc_k$. 
\end{defn}

\begin{defn}
A \emph{LB space} (respectively \emph{LF space}) over $k$ is a locally convex topological vector space over $E$ such that
\begin{equation} \label{eqn:LF}
E \cong \limind_{n \in \Nb} E_n
\end{equation}
where $E_n$ is a Banach space (respectively a  Fr\'echet space) and the system maps are injective. 
\end{defn}

All LB spaces are normal locally convex spaces but the functor ${(-)}^b$ may not commute with the colimit (\ref{eqn:LF}). 

\begin{defn} \label{defn:regular_LF}
	An LB space or an LF space over $k$ is said to be \emph{regular} if 
	\[ E^b \cong \limind_{n \in \Nb} E_n^b. \]
\end{defn}

\begin{defn}
We say that $E \in \ttBorn_k$ is a \emph{bornological LB space} if $E \cong F^b$ where $F \in \ttLoc_k$ is a regular LB space.
\end{defn}

The following criterion will be very important in following sections.

\begin{prop} \label{prop:limpro_commutation_t}
	Consider a projective system $\{ E_{i} \}_{i \in I}$ of objects of $\ttNBorn_{k}$. The following conditions are equivalent
	\begin{enumerate}
		\item $\underset{i \in I}\limpro (E_i^t)$ is normal in $\ttLoc_k$;
		\item $\underset{i \in I}\limpro (E_i^t )\cong (\underset{i \in I}\limpro E_i)^t$.
	\end{enumerate}
\end{prop}

{\bf Proof.}
If $\underset{i \in I}\limpro (E_i^t)$ is normal, then by Observation \ref{obs:tbequiv}
\[ \limpro_{i \in I} (E_i^t) \cong ((\limpro_{i \in I} E_i^t)^b)^t \cong (\limpro_{i \in I} (E_i^{tb}))^t \cong (\limpro_{i \in I} E_i)^t. \]

On the other hand, the same calculations also shows that if $\underset{i \in I}\limpro( E_i^t) \cong (\underset{i \in I}\limpro E_i)^t$
\[ ((\limpro_{i \in I} E_i^t)^b)^t \cong (\limpro_{i \in I} (E_i^{tb}))^t \cong (\limpro_{i \in I} E_i)^t \cong \limpro_{i \in I} (E_i^t) .\]
Therefore, $\underset{i \in I}\limpro (E_i^t)$ is normal.
\ \hfill $\Box$

Finally, in next section we will use the following class of bornological vector spaces.

\begin{defn}\label{defn:regular}
	A bornological vector space is said to be \emph{regular} if it has a basis for its bornology made of subsets which are closed subsets with respect to the topology of $E^t$.
\end{defn}

\begin{rem}
	The condition of Definition \ref{defn:regular} is strictly stronger than the condition of properness of Definition \ref{def:proper_born}, where the requirement is that $E$ has a basis of bornologically closed subsets.
\end{rem}

\subsection{Nuclear bornological spaces and flatness} \label{sec:nuclear}

The functor $\ttCBorn_k \to \ttCBorn_k$, given by $E \mapsto E \wotimes_{\pi, k} F$ where $F$ is a complete bornological vector space, is right exact (even strongly right exact) but not left exact, in general. In order to have a sufficient condition to ensure the exactness of such functors, we use the notion of nuclear operator. This notions was introduced by Grothendieck in the context of locally convex spaces over Archimedean base fields. We will also need to introduce the injective tensor product over maximally complete fields (non-Archimedean or not). 

\begin{defn} \label{defn:BanNucMap}
	Conside a morphism $f: E \to F$ in $\ttBan_{k}$. The morphism $f$ is called \emph{nuclear} if there exist two sequences $\{ \alpha_n \} \subset \uHom(E, k)$ and $\{ f_n \} \subset F$ where
	\begin{itemize}
	\item
$\underset{n = 0}{\overset{\infty}\sum} \|\alpha_n\| \|f_n\| < \infty$
	if $k$ is Archimedean;
	\item
$\|\alpha_n\| \|f_n\| \to 0, \text{ for } n \to \infty $
	if $k$ is non-Archimedean;
\end{itemize}	
such that 
	\[ f(x) = \sum_{n = 0}^\infty \alpha_n(x) f_n  \]
	for all $x \in E$.
\end{defn}

\begin{rem}\label{rem:CatDefCompose}A morphism $f$ is nuclear in the sense of Definition \ref{defn:BanNucMap} if and only if $f$ is in the image of the canonical morphism 
\[ \Hom(k, E^{\vee}\wotimes_{k}F) \to \Hom(k,\uHom(E,F)) \cong \Hom(E, F), \]
where $E^{\vee} = \uHom(E,k)$. This definition makes sense in any closed, symmetric monoidal category. In this abstract context, the pre-composition of a nuclear morphism with any morphism is a nuclear morphism. Similarly, the post-composition of a nuclear morphism with any morphism is nuclear (see Proposition 2.1 of \cite{HR} for a proof of these facts).
\end{rem}

\begin{defn} \label{defn:compactoid}
Let $V$ be a locally convex $k$-vector space. A subset $B \subset V$
is called \emph{compactoid} if for any neighborhood of the origin $U$ there is a finite subset $F \subset V$
such that 
\[ B \subset U + \Gamma(F), \]
where $\Gamma$ denotes the absolute convex hull of $F$.
\end{defn}

Recall that a subset $B \subset V$ of a locally convex $k$-vector space is called \emph{pre-compact} if for any neighborhood of the origin $U$
\[ B \subset U + F \]
for a finite set $F \subset V$.

\begin{prop}\label{prop:CompactoidPrecompact}
Let $V$ be a locally convex $k$-vector space. If $k$ is locally compact then a subset $B \subset V$ is compactoid if and only if it is pre-compact. 
\end{prop}
{\bf Proof.}
Let $B \subset V$ be a compactoid subset and $U \subset V$ be an absolutely convex neighborhood of zero. Then, there is a finite set $F = \{x_1, \ldots, x_n \} \subset V$ such that $B \subset U + \Gamma(F)$. The map $k^n \to E$ given by $(\la_1, \dots, \la_n) \mapsto \underset{i = 1}{\overset{n}\sum} \la_i x_i$ is continuous and it maps the compact set $(k^{\circ})^n$ onto $\Gamma(F)$. Therefore, $\Gamma(F)$ is compact. This means that there is a finite set $F' \subset V$ such that $\Gamma(F) \subset U + F'$, so
\[ B \subset U + U + F'. \]
If $k$ is non-Archimedean we can conclude that 
\[ B \subset U + F' \]
proving that $B$ is pre-compact. If $k$ is Archimedean, we can notice that the same reasoning can be done with $\frac{U}{2}$ in place of $U$, so that
\[ B \subset \frac{U}{2} + \frac{U}{2} + F' \subset U + F' \]
proving that $B$ is pre-compact also in this case.

The converse statement, that every pre-compact subset is compactoid, is obvious.
\ \hfill $\Box$

Thanks to Proposition \ref{prop:CompactoidPrecompact}, a compactoid subset in the case when $k$ is Archimedean is simply pre-compact subset. The notion of compactoid is necessary when the base field is not locally compact, because in this case the definition of pre-compact subset is not useful.

\begin{defn}
The family of compactoid subsets of a locally convex topological vector space $V$ over $k$ forms a bornology, which is denoted by $\ttCpt(V)$.
\end{defn}

Compactoid subsets are always bounded subsets in the sense of von Neumann, but the converse is often false. For example, one can show that if on a Banach space every bounded subset is compactoid, then the Banach space must be finite dimensional.

\begin{defn} \label{defn:compactoid_map}
	Consider a morphism $f: E \to F$ in $\ttBan_{k}$. The morphism $f$ is called \emph{compactoid} if the image of the unit ball of $E$ is a compactoid (in the sense of Definition \ref{defn:compactoid}) subset of $F$.
\end{defn}

\begin{rem}\label{rem:BoundedMor} A morphism $f: E \to F$ in $\ttBan_{k}$ is compactoid if and only if it is a bounded morphism when considered as a map  of bornological spaces $f: E^{b} \to \ttCpt(F)$.
\end{rem}

Notice that what we called a compactoid map in literature is just called a compact map, in the case in which $k$ is Archimedean. If $k$ is Archimedean a nuclear morphism is compact(oid) but there exist compact(oid) morphisms which are not nuclear. However, the non-Archimedean case is simpler.

\begin{prop} \label{prop:compactor_nuclear_map}
Let $k$ be non-Archimedean and $f: E \to F$ a bounded morphism in $\ttBan_{k}$. Then $f$ is nuclear if and only if it is compactoid.
\end{prop}
{\bf Proof.}
In classical terminology what we called nuclear maps over a non-Archimedean field are called completely continuous maps. With this terminology, the proposition is proved as Proposition 2 in page 92 of \cite{Gruson}.
\ \hfill $\Box$

\begin{defn}\label{defn:nuclear_born}
	An object of $\ttCBorn_{k}$ is said to be \emph{nuclear} if there exists an isomorphism
	\[ E \cong \limind_{i \in I} E_i \]
	as in Definition \ref{defn:complete_born} where for any $i < j$ in $I$ the map $E_i \to E_j$ is a nuclear monomorphism of $k$-Banach spaces. Therefore, by definition, nuclear bornological spaces are always complete bornological spaces. 
\end{defn}

We now limit the discussion to the case when $k$ is maximally complete. This restriction will be removed later.

\begin{defn}\label{defn:injective_norm}
Let $k$ be a maximally complete valuation field. Let $E, F$ be objects of $\ttBan_{k}$ and let with $E^{\vee}=\uHom_{\ttBan_{k}}(E,k)$ and $F^{\vee}=\uHom_{\ttBan_{k}}(F,k)$. Then, the \emph{(completed) injective tensor product} is defined to be the completion of the algebraic tensor product $E\otimes_{k}F$ equipped with the semi-norm
	\[  \| \sum e_i \otimes f_i \|_{\epsilon,k} = \underset{\alpha \in (E^{\vee})^\circ, \beta \in (F^{\vee})^\circ}\sup |\sum \alpha(e_i) \beta(f_i)| \]
	where $(E^{\vee})^\circ$ and $(F^{\vee})^\circ$ denote the unit balls. We denote the injective tensor product by $E \otimes_{\epsilon,k} F$ ad the completed one by $E \wotimes_{\epsilon,k} F$.
\end{defn}

The following is Theorem 1 of page 212 of \cite{H2}.
\begin{lem}\label{lem:DifferentTensors}
Let $k$ be a non-Archimedean, maximally complete valuation field. Then, for any $E,F \in \ttBan_{k}$, the natural morphism in $\ttBan_{k}$
	\[ E \otimes_{\pi,k} F \to E \otimes_{\epsilon,k} F \]
	is an isomorphism.
\end{lem}

\begin{defn}\label{defn:ITPborn}
	Let $k$ be a maximally complete valuation field and let $E, F$ be two objects of $\ttCBorn_{k}$. We define the \emph{injective tensor product} of $E$ and $F$ as the colimit in $\ttCBorn_{k}$ of an the  monomorphic inductive system:  
	\begin{equation} \label{eqn:inj_tens_born}
	E \otimes_{\epsilon,k} F = \limind_{B_E \in \mD_E, B_F \in \mD_F} E_{B_E} \otimes_{\epsilon,k} F_{B_F},
	\end{equation}
	where $\mD_E$ and $\mD_F$ are the family of bounded disks of $E$ and $F$ respectively. The \emph{completed injective tensor product} of $E$ and $F$ is defined
	\[ E \wotimes_{\epsilon,k} F = \what{E \otimes_{\epsilon,k} F}. \]
\end{defn}

\begin{rem}
The inductive system (\ref{eqn:inj_tens_born}) is monomorphic because $(-) \otimes_{\epsilon,k} (-)$ is a strongly left exact functor, cf. Corollary 2 page 210 of \cite{H2}.
\end{rem}

Definition \ref{defn:ITPborn} differs from the one given in \cite{H2} chapter 4, \S 2, but is equivalent to it for thanks to the following lemma.

\begin{lem} \label{lemma:limind_injective_tensor}
Let $k$ be a maximally complete valuation field. 
	Let $\{ E_i \}_{i \in I}$ and $\{ F_j \}_{j \in J}$ be two filtered inductive systems of complete bornological vector spaces such that all $E_i$ and $F_i$ are separated and regular (in the sense of Definition \ref{defn:regular}) and all system morphisms are injective. Then, in $\ttCBorn_{k}$ we have that
	\[ E \otimes_{\epsilon,k} F \cong \limind_{(i,j) \in I \times J} ( E_i \otimes_{\epsilon,k} F_i ) \cong ( \limind_{i \in I} E_i ) \otimes_{\epsilon,k} (\limind_{j \in J} F_j), \]
where on the right hand side the injective tensor product is the one considered in \cite{H2}.
\end{lem}
{\bf Proof.}
See \cite{H2} Proposition 4, page 210.
\ \hfill $\Box$

\begin{rem}
The injective tensor product does not commute with filtered inductive limits in general.
\end{rem}

Since all Banach spaces are obviously regular and complete bornological vector spaces are essentially monomorphic objects in $\ttInd(\ttBan_k)$, Lemma \ref{lemma:limind_injective_tensor} readily implies that our definition of injective tensor product agrees with the definition of \cite{H2} page 208, for proper complete bornological vector spaces.

\begin{prop} \label{prop:nuclear_proj_equal_injecive}
	Let $k$ be a maximally complete valuation field.  Let $E$ and $F$ be objects of $\ttCBorn_{k}$. In the Archimedean case we assume that $F$ is nuclear. Then, the natural morphism 
	\[ E \wotimes_{\pi,k} F \longrightarrow E \wotimes_{\epsilon,k} F \]
	in $\ttCBorn_{k}$ is an isomorphism.
\end{prop}

{\bf Proof.}
We start with the case in which $k$ is non-Archimedean.  Write $E=\underset{i \in I}\limind E_i$ and $F=\underset{j \in J}\limind F_j$ as in Definition \ref{defn:complete_born}. The isomorphisms $E_{i} \otimes_{\pi,k} F_{j} \to E_{i} \otimes_{\epsilon,k} F_{j}$ from Lemma \ref{lem:DifferentTensors} gives an isomorphism of systems
\[ \limind_{(i,j) \in I \times J}(E_{i} \otimes_{\pi,k} F_{j}) \to \limind_{(i,j) \in I \times J} (E_{i} \otimes_{\epsilon,k} F_{j}).
\]
The projective tensor product commutes with all colimits and by Lemma \ref{lemma:limind_injective_tensor} also the injective tensor product commutes with colimit we are calculating. Hence, we get an isomorphism 
\[(\limind_{i \in I} E_{i} )\otimes_{\pi,k} (\limind_{j \in J} F_{j} )\to ( \limind_{i \in I} E_{i}) \otimes_{\epsilon,k}  (\limind_{j \in J} F_{j}) 
\]
and the functoriality of the completion yields the claimed isomorphism.

We now consider the case in which $k$ is Archimedean. First recall the following fact: let $f: E_1 \to E_2$ be a nuclear morphism of Banach spaces and $V$ another Banach space, then the linear map of vector spaces
\[f \otimes_{k} \text{id}_V: E_1 \otimes_{\epsilon,k} V \to E_2 \otimes_{\pi,k} V \] 
is bounded. This is precisely the content of Proposition 2.4, page I-15 of \cite{DV}, where the complex case is discussed, but the same arguments work also over $\Rb$. Then, consider $F = \limind\limits_{j \in J} F_j$ and  $E = \limind\limits_{i \in I} E_i$ presentations of $F$  and $E$ as in Definition \ref{defn:complete_born} and as in Definition \ref{defn:nuclear_born} respectively. For any $(i_1, j_1) < (i_2, j_2)$ the morphisms
\[ \a_{i_1, j_1}: E_{i_1} \otimes_{\epsilon,k} F_{j_1} \to E_{i_2} \otimes_{\epsilon,k} F_{j_1} \to E_{i_2} \otimes_{\pi,k} F_{j_2} \]
are bounded morphisms of Banach spaces because of the mentioned Proposition 2.4 of \cite{DV}. By Lemma \ref{lemma:limind_injective_tensor} we have that the map
\[\begin{split} E \otimes_{\epsilon,k} F &\cong (\limind_{i \in I} E_i) \otimes_{\epsilon,k} (\limind_{j \in J} F_j) \cong \limind_{(i,j) \in I \times J} (E_i \otimes_{\epsilon,k} F_j)  \stackrel{\a}{\to} \\  & \limind_{(i,j) \in I \times J} (E_i \otimes_{\pi,k} F_j) \cong (\limind_{i \in I} E_i) \otimes_{\pi,k} (\limind_{j \in J} F_j) \cong E \otimes_{\pi,k} F
\end{split} \]
obtained by the morphisms $\a_{i_1, j_1}$ is a bounded map, and is easy to check that the composition of $\a$ with the canonical map
\[ E \otimes_{\pi,k} F \to E \otimes_{\epsilon,k} F \]
is the identity on both sides. We conlcude the proof using again the functorialiy of the completion of bornological vector spaces.
\ \hfill $\Box$

\begin{lem} \label{lem:injective_tensor_strict_mono} 
Let $k$ be a maximally complete field and $F \in \ttCBorn_k$. The functor $(-) \otimes_{\epsilon,k} F$ is strongly left exact.
\end{lem}
{\bf Proof.}
See \cite{H2} Corollary 2, page 210.
\ \hfill $\Box$

\begin{lem} \label{lem:exact_completion} 
The completion functor $\widehat{(-)}: \ttNrm_k \to \ttBan_k$ is exact.
\end{lem}
{\bf Proof.}
Proposition 4.1.13 of \cite{Pr2} shows that the completion functor is exact for locally convex spaces. Since the functor $\widehat{(-)}: \ttNrm_k \to \ttBan_k$ precisely agrees with the completion as locally convex spaces, we can use that proposition to deduce our lemma. Notice that \cite{Pr2} deals only with Archimedean base fields but the same reasoning works also for non-Archimedean base fields since only the uniform structures of the spaces are considered.
\ \hfill $\Box$

\begin{rem}
$\widehat{(-)}: \ttNrm_k \to \ttBan_k$ is not strongly left exact and does not preserve monomorphisms.
\end{rem}

\begin{lem}\label{lem:nArchBanFlat} 
If $k$ is non-Archimedean, then any object in the category of non-Archimedean Banach spaces over $k$ is flat.
\end{lem}
{\bf Proof.} 
If $k$ is maximally complete, then the claim follows from the combination of Lemma \ref{lem:injective_tensor_strict_mono}, Lemma \ref{lem:DifferentTensors} and Lemma \ref{lem:exact_completion} .

If $k$ is not maximally complete we can choose a maximal completion $K/k$. For every $E \in \ttBan_k$ the canonical map $\iota_E: E \to E \otimes_k K$  is a strict monomorphism (even an isometry onto its image, cf. Lemma 3.1 of \cite{Poi} applied with $A = k, B = K$ and $C = E$).
Consider a strict exact sequence
\begin{equation} \label{eqn:flat_ban}
 0 \to \ker f \to E \stackrel{f}{\to} F
\end{equation}
in $\ttBan_k$. We can always suppose that $\ker f$ is equipped with the restriction of the norm of $E$, in which case $\ker f \to E$ is an isometric embedding. Applying again Lemma 3.1 of \cite{Poi}, we obtain a diagram
\[
\begin{tikzpicture}
\matrix(m)[matrix of math nodes,
row sep=2.6em, column sep=2.8em,
text height=1.5ex, text depth=0.25ex]
{  0  & \ker f & E & F  \\
   0  & \ker f \otimes_k K & E \otimes_k K & F \otimes_k K  \\};
\path[->,font=\scriptsize]
(m-1-1) edge node[auto] {$$} (m-1-2);
\path[->,font=\scriptsize]
(m-2-1) edge node[auto] {$$} (m-2-2);
\path[->,font=\scriptsize]
(m-1-2) edge node[auto] {$$} (m-1-3);
\path[->,font=\scriptsize]
(m-1-2) edge node[auto] {$$} (m-2-2);
\path[->,font=\scriptsize]
(m-2-2) edge node[auto] {$$}  (m-2-3);
\path[->,font=\scriptsize]
(m-1-3) edge node[auto] {$$}  (m-2-3);
\path[->,font=\scriptsize]
(m-1-3) edge node[auto] {$f$}  (m-1-4);
\path[->,font=\scriptsize]
(m-2-3) edge node[auto] {$$}  (m-2-4);
\path[->,font=\scriptsize]
(m-1-4) edge node[auto] {$$}  (m-2-4);
\end{tikzpicture}
\]
where all vertical maps are isometric embeddings. The bottom row is exact algebraically. It is also strictly exact, because we can apply  Lemma 3.1 of \cite{Poi} to the isometric embedding $\ker f \to E$ (precisely using $A = \ker f, B = E$ and $C = K$ this time) to deduce that $\ker f \otimes_k K \to E \otimes_k K$ is an isometric embedding. Then, consider and a $G \in \ttBan_k$. Tensoring (\ref{eqn:flat_ban}) with $G$ we obtain a diagram
\[
\begin{tikzpicture}
\matrix(m)[matrix of math nodes,
row sep=2.6em, column sep=2.8em,
text height=1.5ex, text depth=0.25ex]
{  0  & \ker f \otimes_k G & E \otimes_k G & F \otimes_k G  \\
	0  & (\ker f \otimes_k G) \otimes_k K & (E \otimes_k G) \otimes_k K & (F \otimes_k G) \otimes_k K  \\};
\path[->,font=\scriptsize]
(m-1-1) edge node[auto] {$$} (m-1-2);
\path[->,font=\scriptsize]
(m-2-1) edge node[auto] {$$} (m-2-2);
\path[->,font=\scriptsize]
(m-1-2) edge node[auto] {$$} (m-1-3);
\path[->,font=\scriptsize]
(m-1-2) edge node[auto] {$$} (m-2-2);
\path[->,font=\scriptsize]
(m-2-2) edge node[auto] {$$}  (m-2-3);
\path[->,font=\scriptsize]
(m-1-3) edge node[auto] {$$}  (m-2-3);
\path[->,font=\scriptsize]
(m-1-3) edge node[auto] {$f$}  (m-1-4);
\path[->,font=\scriptsize]
(m-2-3) edge node[auto] {$$}  (m-2-4);
\path[->,font=\scriptsize]
(m-1-4) edge node[auto] {$$}  (m-2-4);
\end{tikzpicture}
\]
where the top row is algebraically exact. We need only to check that $\ker f \otimes_k G \to E \otimes_k G$ is a strict monomorphism.
Using the isomorphisms
\begin{eqnarray*} (E \otimes_k G) \otimes_k K \cong (E \otimes_k K) \otimes_K (G \otimes_k K), \\ (F \otimes_k G) \otimes_k K \cong (F \otimes_k K) \otimes_K (G \otimes_k K) \\
(\ker f \otimes_k G) \otimes_k K \cong (\ker f \otimes_k K) \otimes_K (G \otimes_k K) 
\end{eqnarray*} 
last diagram becomes 
\[
\begin{tikzpicture}
\matrix(m)[matrix of math nodes,
row sep=2.6em, column sep=2.8em,
text height=1.5ex, text depth=0.25ex]
{  0  & \ker f \otimes_k G & E \otimes_k G & F \otimes_k G  \\
	0  & (\ker f \otimes_k K) \otimes_K (G \otimes_k K) & (E \otimes_k K) \otimes_K (G \otimes_k K) & (F \otimes_k K) \otimes_K (G \otimes_k K)  \\};
\path[->,font=\scriptsize]
(m-1-1) edge node[auto] {$$} (m-1-2);
\path[->,font=\scriptsize]
(m-2-1) edge node[auto] {$$} (m-2-2);
\path[->,font=\scriptsize]
(m-1-2) edge node[auto] {$$} (m-1-3);
\path[->,font=\scriptsize]
(m-1-2) edge node[auto] {$$} (m-2-2);
\path[->,font=\scriptsize]
(m-2-2) edge node[auto] {$$}  (m-2-3);
\path[->,font=\scriptsize]
(m-1-3) edge node[auto] {$$}  (m-2-3);
\path[->,font=\scriptsize]
(m-1-3) edge node[auto] {$f$}  (m-1-4);
\path[->,font=\scriptsize]
(m-2-3) edge node[auto] {$$}  (m-2-4);
\path[->,font=\scriptsize]
(m-1-4) edge node[auto] {$$}  (m-2-4);
\end{tikzpicture}
\]
where the bottom row is strictly exact because $K$ is maximally complete. Hence, in the diagram
\[
\begin{tikzpicture}
\matrix(m)[matrix of math nodes,
row sep=2.6em, column sep=2.8em,
text height=1.5ex, text depth=0.25ex]
{   \ker f \otimes_k G & E \otimes_k G  \\
	 (\ker f \otimes_k K) \otimes_K (G \otimes_k K) & (E \otimes_k K) \otimes_K (G \otimes_k K)  \\};
\path[->,font=\scriptsize]
(m-1-1) edge node[auto] {$$} (m-1-2);
\path[->,font=\scriptsize]
(m-2-1) edge node[auto] {$$} (m-2-2);
\path[->,font=\scriptsize]
(m-1-1) edge node[auto] {$$} (m-2-1);
\path[->,font=\scriptsize]
(m-1-2) edge node[auto] {$$} (m-2-2);
\end{tikzpicture}
\]
all maps are known to be strict monomorphism but the top horizontal, which is therefore a strict monomorphism too. This shows that the functor $(-) \otimes_{\pi, k} (-)$ is strongly exact, and applying Lemma \ref{lem:exact_completion}  we deduce that $(-) \wotimes_{\pi, k} (-)$ is exact.
\ \hfill $\Box$

We now look at flatness in the closed symmetric monoidal category $\ttCBorn_{k}$.
\begin{thm} \label{thm:strct_exact_nuclear}
	Let $F$ be complete bornological vector spaces over $k$. If $k$ is non-Archimedean then $F$ is flat with respect the category $(\ttCBorn_{k}, \wotimes_{k},k)$. If $k$ is Archimedean then the same conclusion holds provided that $F$ is nuclear.
\end{thm}

{\bf Proof.}
When $k$ is Archimedean, it is of course maximally complete. So, the theorem follows immediately from combining Lemma \ref{lem:injective_tensor_strict_mono} and Proposition \ref{prop:nuclear_proj_equal_injecive}. 

Suppose now that $k$ is a non-Archimedean complete valuation field and fix a representation $F = \underset{j \in J}\limind F_j$ as a monomorphism filtered inductive system of Banach spaces. Unraveling the definitions, the functor $(-) \wotimes_{\pi, k} F$ can be written as
\[ (-) \wotimes_{\pi, k} F = 
((-) \otimes_{\pi, k} F)^{\wedge} = (\limind_{j \in J}(-) \otimes_{\pi, k} F_j)^{\wedge} = \limind_{j \in J}(((-) \otimes_{\pi, k} F_j)^{\wedge})\]
where last colimit is calculated in $\ttCBorn_k$, meaning that it is the separated colimit of bornological vector spaces. Since $\ttCBorn_k$ is an elementary quasi-abelian category, the filtered colimits are exact. Therefore, applying Lemma \ref{lem:nArchBanFlat} and Lemma \ref{lem:exact_completion} we see that the functor $(-) \wotimes_{\pi, k} F$ can be written as a composition of exact functors.
\ \hfill $\Box$

We conclude this section with some results about nuclearity needed in subsequent sections.

\begin{prop} \label{prop:nuclear_proj}
	Any countable projective limit in $\ttCBorn_{k}$ of nuclear objects is nuclear (and therefore by Theorem \ref{thm:strct_exact_nuclear} flat).
\end{prop}
{\bf Proof.}
The proposition is proved in \cite{Hog3} Theorem 4.1.1, page 200, for the case in which $k$ is an Archimedean base field. So, we check the claim only for non-Archimedean base fields.

First we show that a closed subspace of a nuclear bornological space is nuclear. Let $E \subset F$ be a closed subspace of a nuclear space and let $F \cong \underset{i \in I}\limind F_i$ as in Definition \ref{defn:nuclear_born}. As bornological vector spaces, there is an isomorphism 
\[ E \cong \limind_{i \in I} ( E \cap F_i ) \]
and $E \cap F_i$ are $k$-Banach spaces. It is enough to check that the morphisms $E \cap F_i \to E \cap F_j$ are nuclear for each $i < j$. Thanks to Proposition \ref{prop:compactor_nuclear_map} this is equivalent to check that the maps $E \cap F_i \to E \cap F_j$ are compactoid. This follows from \cite{PGS}, Theorem 8.1.3 (ii) and (iii).

Now, we check that countable products of nuclear bornological spaces are nuclear. Thus, suppose that $\{ E_n \}_{n \in \Nb}$ is a family of nuclear bornological spaces and put $E = \underset{n \in \Nb}\prod E_n$. Fix for any $n$ a completant bounded disk $B_n \subset E_n$ such that there exists a completant bounded disk $A_n \subset E_n$ for which $B_n \subset A_n$ and the inclusion $E_{B_n} \to E_{A_n}$ is a nuclear map (or equivalently compactoid, cf. Proposition \ref{prop:compactor_nuclear_map}). Let $B = \underset{n \in \Nb}\prod B_n$ and $A = \underset{n \in \Nb}\prod A_n$. It is clear that by varying $B_n$ over a final family of bounded completant disks of $E_n$, for each $n$, the family of disks $B \subset E$ and $A \subset E$ obtained in this way form a final family of bounded disks for the bornology of $E$. By construction the map $E_B \to E_A$ is bounded, where $E_A$ is the subspace of $\underset{n \in \Nb}\prod E_{A_n}$ consisting of bounded sequences equipped with the supremum norm, and $E_B$ is the analogous subspace of $\underset{n \in \Nb}\prod E_{B_n}$. We need to check that $E_B \to E_A$ sends the unit ball of $E_B$ to a compactoid subset of $E_A$. Fix a $\a \in |k^\times|$ such that $\a < 1$. By rescaling the norms of $E_{B_n}$ and $E_{A_n}$ we can assume that $B_n$ is contained in the ball of radius $\a^n$ of $E_{A_n}$.
Since $B_n$ is compactoid in $E_{A_n}$, by Theorem 3.8.25 of \cite{PGS} there exists a zero sequence $\{ x_i^{(n)} \}_{i \in \Nb}$ such that $B_n \subset \Gamma(\{ x_i^{(n)} \}_{i \in \Nb})$, where $\Gamma$ denotes the absolutely convex hull. Then, we can consider the sequence
\[ x^{(n)} = \{ (x_i^{(n)}) \}_{i \in \Nb} \subset E_A. \]
Since by hypothesis, for each $i$ we have that
\[ \lim_{n \to \infty} |x_i^{(n)}|_{E_{A_n}} = 0 \]
then the sequence $x^{(n)}$ is a zero sequence of $E_A$ whose absolutely convex hull is $B$. So applying Theorem 3.8.24 of \cite{PGS} we can deduce that $B$ is compactoid in $E_A$.
\ \hfill $\Box$ 

\begin{prop} \label{prop:nuclear_ind_lim} 
	Any small inductive limit in $\ttCBorn_{k}$ of nuclear objects is nuclear.
\end{prop}

{\bf Proof.}
It is an easy consequence of Definition \ref{defn:nuclear_born} that any monomorphic filtered inductive limit of bornological nuclear spaces is a nuclear bornological space. Then, by the description of coproducts in $\ttCBorn_k$ of Lemma 2.7 of \cite{BaBe} (which agree with coproducts in $\ttInd(\ttBan_k)$) and applying Proposition \ref{prop:nuclear_proj} to finite coproducts, we obtain that coproducts of bornological nuclear spaces are nuclear bornological spaces. Therefore, it remains only to check that quotients of nuclear spaces are nuclear.

The Archimedean case of the proposition is proved in Theorem 4.1.1, page 200 of \cite{Hog3}. So, let $k$ be non-Archimedean, $F$ a nuclear bornological space over $k$ and $E \subset F$ a closed subspace. Fixing an isomorphism $F \cong \underset{i \in I}\limind F_i$ as in Definition \ref{defn:nuclear_born}, we can write
\[ \frac{F}{E} \cong \limind_{i \in I} \frac{F_i}{F_i \cap E}. \]
Thanks to Proposition \ref{prop:compactor_nuclear_map} it is enough to prove that the system morphisms $\phi_{i, j}: \frac{F_i}{F_i \cap E} \to \frac{F_j}{F_j \cap E}$ for any $i \le j$ are compactoid maps. By hypothesis $F_i \to F_j$ is a compatoid map and since the image of compactoid subsets by bounded maps are compactoid subsets, then also $F_i \to \frac{F_j}{F_j \cap E}$ is a compactoid map. By Theorem 8.1.3 (xi) of \cite{PGS} this is equivalent to say that $\phi_{i,j}$ is a compactoid map.
\ \hfill $\Box$

We recall that a dagger affinoid algebra is a bornological algebra which is isomorphic to a quotient of the algebra $\mW_k^n(\rho)$ of over-convergent analytic functions on the polycylinder of polyradius $\rho$. The isomorphism
\[ \mW_k^n(\rho) \cong \limind_{r > \rho} \mT_k^n(r) \]
where $\mT_k^n(r)$ are the $k$-Banach algebras of strictly convergent power-series on the polycylinder of polyradius $\rho$, endow $\mW_k^n(\rho)$ with a bornology and dagger affinoid algebras are endowed with the quotient bornology. We refer to Section4 of \cite{BaBe} of chapter 3 of \cite{Bam} for more details.

\begin{prop} \label{prop:dagger_nuclear} 
The underlying bornological vector spaces of dagger affinoid algebras are nuclear.
\end{prop}

{\bf Proof.}
Since by Proposition \ref{prop:nuclear_ind_lim} quotients of bornological nuclear spaces are nuclear, it is enough to check that $\mW_k^n(\rho)$ is nuclear. We first look at the case in which $k$ is non-Archimedean. It is easy to check the canonical restriction maps $\mT_k^n(\rho) \to \mT_k^n(\rho')$, for $\rho' < \rho$ are compactoid maps, which proves the claim. For details the reader can see Theorem 11.4.2 and Remark 11.4.3 of \cite{PGS} where there $\mA^\dagger(\rho)$ is what in our notation is $\mW_k^n(\rho)$.

If $k$ is Archimedean one can write $\mW_k^n(\rho)$ as the direct limit of the filtered system of Fr\'echet spaces of holomorphic functions on open polydisks of polyradius bigger than $\rho$, see Section3.3 of \cite{Bam} for a detailed proof of this fact. Thanks to a celebrated theorem of Montel these Fr\'echet spaces are nuclear Fr\'echet spaces (see for example Proposition 4 on page 26 of \cite{DV} for a proof of this fact) and by Lemma \ref{lem:frechet_nuclear} this is equivalent to say that these spaces are nuclear as bornological spaces. So, since by Lemma \ref{prop:nuclear_ind_lim} direct limits of nuclear bornological spaces are nuclear we deduce that $\mW_k^n(\rho)$ is nuclear.
\ \hfill $\Box$

\subsection{A relative flatness lemma} \label{sec:relative_flat}

In last section we proved that if $F$ is a nuclear bornological vector space then the functor $(-) \wotimes_{k} F$ is exact. Now we want to find conditions on $(-) \wotimes_{k} F$ to preserve a different kind of projective limit: The cofiltered ones. The notion of projective tensor product was introduced by Grothendieck, in the category of locally convex spaces, in order to find a topological tensor product which commutes with projective limits. This is the origin of the name projective tensor product. But the bornological projective tensor product does not commute with projective limits, in general. So, we need to find some sufficient conditions to ensure this commutation for the projective limits we will study in Section\ref{sec:Stein_geometry}. This study involves comparing the bornological and topological projective tensor products. Thus, we start by recalling the notion of topological projective tensor product.

\begin{defn} \label{defn:top_proj_tensor}
Given $E, F \in \ttLoc_k$ one defines $E \otimes_{\pi,k} F \in \ttLoc_k$ as the algebraic tensor product equipped with the locally convex topology whose base of neighborhoods of zero is given by the family of absolutely convex hulls of subsets of $E \otimes_{k} F$ of the form
\[ U \otimes V = \{ x \otimes y \ | \ x \in U, y \in V \} \]
where $U$ and $V$ vary over the family of neighborhoods of $E$ and $F$ respectively. The \emph{complete projective tensor product} is defined as the separated completion of the space $E \otimes_{\pi,k} F$. The space obtained in this way is denoted $E \wotimes_{\pi,k} F$.
\end{defn}

\begin{defn} \label{defn:top_nuclear}
$E \in \ttLoc_k$ is said to be \emph{nuclear} if $E \cong \underset{i \in I} \limpro E_i$ for a cofiltered epimorphic projective system of semi-normed spaces $\{E_{i}\}_{i \in I}$ such that, for each $i < j$ in $I$, the maps $\what{E}_j \to \what{E}_i$ of Banach spaces induced on the separated completions are nuclear maps.
\end{defn}

\begin{defn} \label{defn:top_inj_tensor}
	Given $E, F \in \ttLoc_k$ one defines  the \emph{injective tensor product} $E \otimes_{\epsilon,k} F \in \ttLoc_k$ as the algebraic tensor product equipped with the semi-norms $\{ \| \cdot \|_{p_i \otimes_\epsilon q_j} \}_{i,j \in I \times J}$, where for each $i,j \in I \times J$ the semi-norm $\| \cdot \|_{p_i \otimes_\epsilon q_j}$ is defined as in Definition \ref{defn:injective_norm}, where $\{ p_i \}_{i \in I}$ is a family of semi-norms that defines the topology of $E$ and $\{ q_j \}_{j \in J}$ is a family of semi-norms that defines the topology of $F$. The separated completion of $E \otimes_{\epsilon,k} F$ is by definition the \emph{complete injective tensor product} and is denoted $E \wotimes_{\epsilon,k} F$.
\end{defn}

\begin{lem} \label{lem:top_inj_pro}
	If $E \in \ttLoc_k$ is nuclear then the functors $(-) \wotimes_{\epsilon,k} E$ and $(-) \wotimes_{\pi,k} E$ are naturally isomorphic.
\end{lem}
{\bf Proof.}
For the Archimedean case one can refer to Theorem 1 of page 25 of \cite{DV}.
For the non-Archimedean case the result is a direct consequence of Theorem 8.5.1 and Theorem 10.2.7 of \cite{PGS}.
\ \hfill $\Box$

\begin{defn} \label{defn:binuclear}
	Given $E \in \ttNBorn_k$, we say that $E$ is \emph{binuclear} is $E$ is a nuclear bornological space and $E^t$ is a nuclear topological vector space.
\end{defn}

\begin{lem} \label{lem:LB_nuclear}
	Let $E \in \ttLoc_k$. Suppose there exists a countable  monomorphic compactoid inductive system of Banach spaces $\{E_n\}$ be such that $E \cong \underset{n \in \Nb} \limind E_n$  (in particular $E$ is a normal space, cf. Definition \ref{defn:normal_space}). Then, $E$ is nuclear as locally convex space and $E^b$ is nuclear as bornological vector space.
\end{lem}
{\bf Proof.}
If $k$ is Archimedean, then this lemma is proved in \cite{Hog3} Theorem 7, page 160, since compactoid inductive limits are complete. Let $k$ be non-Archimedean. Then, by Theorem 11.3.5 (v) of \cite{PGS} $E$ is a regular LB-space (in the sense of Definition \ref{defn:regular_LF}), hence $E^b \cong \underset{n \in \Nb} \limind E_n$ as bornological vector space, so $E^b$ is nuclear. And by Theorem 11.3.5 (ix) of \cite{PGS} $E$ is nuclear as locally convex space.
\ \hfill $\Box$

For Fr\'echet spaces the situation is a bit more complicated.

\begin{lem} \label{lem:frechet_nuclear}
	If $k$ is Archimedean, then a Fr\'echet space $E \in \ttLoc_k$ is nuclear as locally convex  vector space if and only if $E^b$ is nuclear as bornological vector space. For non-Archimedean base fields, if a Fr\'echet space $E \in \ttLoc_k$ is nuclear as locally convex vector space then $E^b$ is nuclear as bornological space.
\end{lem}
{\bf Proof.}
If $k$ is Archimedean then this lemma is proved in \cite{Hog3} Theorem 7 (i) and (ii), page 160. If $k$ is non-Archimedean then the nuclear Fr\'echet space $E$ is Montel by Corollary 8.5.3 of \cite{PGS}. Therefore, by Theorem 8.4.5 ($\delta$) of \cite{PGS}, each bounded subset of $E$ is compactoid, \ie there is an isomorphism of bornological vector spaces $\ttCpt(E) \cong E^b$. Let's write
\[ E^b \cong \limind_{B \in \mD_E} E_B \]
as in Remark \ref{rem:gauge}, where $\mD_E$ is the family of bounded disks of $E$, and 
\[ \ttCpt(E) \cong \limind_{B \in \ttCpt_E} E_B \]
where $\ttCpt_E$ is the family of compactoid disks of $E$. The family $\ttCpt_E$ is characterized as the family of bounded disks $B$ of $E$ for which there exists another bounded disk $B'$ such that $B \subset B'$ and the canonical map
\[ E_B \to E_{B'} \]
is compactoid. This characterization, in combination with the isomorphism
\[ \limind_{B \in \ttCpt_E} E_B \cong \limind_{B \in \mD_E} E_B \]
imply that for every bounded disk $B$ of $E$ there exists another bounded disk $B'$ such that $B \subset B'$ and the canonical map
\[ E_B \to E_{B'} \]
is compactoid. Since $E$ is complete we can always suppose that $B$ and $B'$ are Banach disks (\ie that $E_B$ and $E_{B'}$ are Banach). This proves that $E^b$ satisfies the conditions of Definition \ref{defn:nuclear_born}.
\ \hfill $\Box$

\begin{rem}
If $k$ is non-Archimedean there exists a Fr\'echet space $E$ for which $E^b$ is nuclear whereas $E$ is not nuclear. This is due to the fact that, for non-Archimedean base fields, a morphism of Banach spaces is nuclear if and only if it is compactoid. But the Archimedean analogous of the notion of compactoid map is what is called a compact map (or operator), in classical functional analysis over $\Rb$ and $\Cb$ . Therefore, a nuclear Fr\'echet space over a non-Archimedean base field, following this terminology, is analogous to what in Archimedean functional analysis is a called Schwartz-Fr\'echet space, (see \cite{TER} and the first chapter of \cite{Hog3} for a detailed account of the properties of Schwartz-Fr\'echet spaces). In the Archimedean case it is known that a Fr\'echet space is bornologically Schwartz if and only if it is Montel (cf. Theorem 8 of page 20 of \cite{Hog3}). Also in the non-Archimedean case a  Fr\'echet space is bornologically nuclear if and only if it is Montel (cf. Corollary 8.5.3 of \cite{PGS}). One can prove (in both settings) that there exist Fr\'echet-Montel spaces which are not Schwartz. For explicit examples of such spaces, see Counterexamples 9.8.2 (vi) of \cite{PGS}, for $k$ non-Archimedean, and \S 4 of \cite{TER} for $k$ Archimedean.
\end{rem}

We will be interested in finding conditions for which the functors ${(-)}^t$ and ${(-)}^b$ intertwine the complete bornological and the complete topological projective tensor products. In general ${(-)}^t$ and ${(-)}^b$ preserve neither the complete nor the incomplete projective tensor products.

\begin{prop} \label{prop:tensor_t_b_frechet}
	Let $E, F \in \ttBorn_k$ be bornological Fr\'echet spaces one of which is binuclear, then
	\[ (E \wotimes_{\pi,k} F)^t \cong E^t \wotimes_{\pi,k} F^t. \]
	If $E, F \in \ttLoc_k$ are Fr\'echet spaces one of which is nuclear, then
	\[ (E \wotimes_{\pi,k} F)^b \cong E^b \wotimes_{\pi,k} F^b. \]
\end{prop}

{\bf Proof.}
In page 215 of \cite{H2} it is proved that for metrizable topological or bornological spaces 
\[ (E \otimes_{\epsilon,k} F)^t \cong E^t \otimes_{\epsilon,k} F^t \]
and 
\[ (E \otimes_{\epsilon,k} F)^b \cong E^b \otimes_{\epsilon,k} F^b. \]
Since for metric topological or bornological spaces the bornological and the topological notion of convergence agree (see the last lines of page 108 of \cite{H2} or Proposition 1.17 of \cite{Bam2} for a detailed proof of this fact) it follows that for metric spaces the notion of bornological and topological completeness agree (see also Corollary 1.18 of \cite{Bam2}). Therefore, we deduce that 
\[ (E \wotimes_{\epsilon,k} F)^t \cong E^t \wotimes_{\epsilon,k} F^t \]
and 
\[ (E \wotimes_{\epsilon,k} F)^b \cong E^b \wotimes_{\epsilon,k} F^b. \]
Finally, since one of $E$ or $F$ is binuclear we deduce that the complete injective and projective tensor products coincide by Lemma \ref{lem:top_inj_pro} and Proposition \ref{prop:nuclear_proj_equal_injecive} (in both categories $\ttLoc_k$ and $\ttBorn_k$), obtaining the required isomorphisms. 
\ \hfill $\Box$

We also underline that in Theorem 2.3 of \cite{M2} the Archimedean version of Proposition \ref{prop:tensor_t_b_frechet} is discussed.

\begin{lem} \label{lem:reduced_loc_conv}
	Let $\{ E_i \}_{i\in I}$ be a cofiltered projective system in $\ttLoc_{k}$ whose projective limit we call $E$ and $F \in \ttLoc_{k}$. Then
	\[ E \wotimes_{\pi,k} F = (\limpro_{i\in I} E_i) 
	\wotimes_{\pi, k} F \cong \limpro_{i \in I}(E_i \wotimes_{\pi, k} F) \]
\end{lem}
{\bf Proof.}
Cf. Proposition 9 at page 192 of \cite{H2}.
\ \hfill $\Box$

\begin{cor} \label{cor:proj_lim_1}
	Let $\{ E_i\}_{i\in \Nb}$ be a cofiltered projective system of Banach vector spaces such that $\underset{i\in \Nb}\limpro E_i$ is a binuclear Fr\'echet bornological vector space and $F$ a Fr\'echet bornological vector space, over $k$. Then, the canonical map
	\[ E \wotimes_{\pi,k} F = (\limpro_{i\in \Nb} E_i) 
	\wotimes_{\pi, k} F \to \limpro_{i \in \Nb}(E_i \wotimes_{\pi, k} F)  \] 
	is an isomorphism of bornological vector spaces.
\end{cor}
{\bf Proof.}
The claim follows form the chian of isomorphisms
\[ E \wotimes_{\pi, k} F \stackrel{\ref{example:normal}}{\cong} (E \wotimes_{\pi, k} F)^{t b} \stackrel{\ref{prop:tensor_t_b_frechet}}{\cong} (E^t \wotimes_{\pi, k} F^t)^b \stackrel{\ref{lem:reduced_loc_conv}}{\cong} (\limpro_{i\in \Nb} (E_i^t \wotimes_{\pi, k} F^t))^b \stackrel{}{\cong} \]
\[  \stackrel{\ref{obs:tbequiv}}{\cong} \limpro_{i\in \Nb} (E_i^t \wotimes_{\pi, k} F^t)^b \stackrel{}{\cong} \limpro_{i\in \Nb} (E_i \wotimes_{\pi, k} F) \]
where last isomorphism is immediate from the definition of projective tensor product (compare Definition \ref{defn:top_proj_tensor} with Definition 3.57 of \cite{BaBe}). 
\ \hfill $\Box$

\begin{cor} \label{cor:proj_lim_2}
Let $\{ E_i\}_{i\in \Nb}$ be a cofiltered projective system of binuclear Fr\'echet bornological vector spaces and $F$ bornological Fr\'echet vector space, over $k$. Then, the canonical map
\[ E \wotimes_{\pi,k} F = (\limpro_{i\in \Nb} E_i) 
\wotimes_{\pi, k} F \to \limpro_{i \in \Nb}(E_i \wotimes_{\pi, k} F)  \] 
is an isomorphism of bornological vector spaces.
\end{cor}
{\bf Proof.}
Notice that Lemma \ref{prop:nuclear_proj}, together with the well-known fact that a projective limit of nuclear spaces in $\ttLoc_k$ is a nuclear space and Proposition \ref{prop:limpro_commutation_t} imply that $E$ is binuclear. So, also for this corollary we have the chian of isomorphisms
\[ E \wotimes_{\pi, k} F \stackrel{\ref{example:normal}}{\cong} (E \wotimes_{\pi, k} F)^{t b} \stackrel{\ref{prop:tensor_t_b_frechet}}{\cong} (E^t \wotimes_{\pi, k} F^t)^b \stackrel{\ref{lem:reduced_loc_conv}}{\cong} (\limpro_{i\in \Nb} (E_i^t \wotimes_{\pi, k} F^t))^b \stackrel{}{\cong} \]
\[  \stackrel{\ref{obs:tbequiv}}{\cong} \limpro_{i\in \Nb} (E_i^t \wotimes_{\pi, k} F^t)^b \stackrel{\ref{prop:tensor_t_b_frechet}}{\cong} \limpro_{i\in \Nb} (E_i \wotimes_{\pi, k} F). \]
\ \hfill $\Box$

 \subsection{Strict exact sequences in bornological and topological settings} \label{sec:exact_sequences}

This section contains some results about how the notions of strictly short exact sequence in $\ttBorn_k$ and $\ttLoc_k$ are related, which are needed later.
 
\begin{lem} \label{lem:compactoid_quotient}
Let $f: E \to F$ be a surjective continuous map between Fr\'echet spaces.  For any compactoid subset $B \subset F$, there exists a compactoid subset $B' \subset E$ such that $f(B') = B$. In particular, the functor $\ttCpt: \ttLoc_k \to \ttCBorn_k$ preserves strict short exact sequences of Fr\'echet spaces.
\end{lem}
{\bf Proof.}
When $k = \Rb, \Cb$, see \cite{M} Theorem 1.62, and for the non-Archimedean case see Theorem 3.8.33 \cite{PGS}.
\ \hfill $\Box$

\begin{lem} \label{lem:compactoid_neumann}
Let $F$ be a nuclear Fr\'echet space, then the von Neumann and the compactoid bornology conicides, \ie the identity map gives an isomorphism $F^b \cong \ttCpt(F)$ in $\ttCBorn_{k}$.
\end{lem}

{\bf Proof.}
In both cases one can show that nuclear Fr\'echet spaces are Montel spaces, and for Montel spaces the von Neumann bornology and the compactoid bornology agree (by definition). Possibile references for these results are sections 8.4 and 8.5 of \cite{PGS} for non-Archimedean base fields and Theorem 8 of page 20 of \cite{Hog3} for Archimedean base fields.
\ \hfill $\Box$

\begin{lem} \label{lem:compactoid_conv}
Let $E$ be a nuclear Fr\'echet space and $\{ x_n \}$ a sequence of elements of $E$. Then,
 the following are equivalent
\begin{enumerate}
\item $\{ x_n \}$ converges topologically to $0$;
\item $\{ x_n \}$ converges bornologically to $0$ for the bornology of $E^b$;
\item $\{ x_n \}$ converges bornologically to $0$ for the bornology of $\ttCpt(E)$.
\end{enumerate}
\end{lem}

{\bf Proof.}
This result is proven for Archimedean base fields in \cite{M2} Corollary 3.8. For the general case: For the equivalence of (1) and (2) we can refer to \cite{Bam2}, Proposition 1.17. Finally, conditions (2) and (3) are equivalent by Lemma \ref{lem:compactoid_neumann}.
\ \hfill $\Box$

\begin{rem}
	We conjecture that Lemma \ref{lem:compactoid_conv} holds for all Fr\'echet spaces, without the nuclearity hypothesis. Indeed the proof given in \cite{M2} for Archimedean base fields, does not use the nuclearity hypothesis but it easily adapts to the non-Archimedean case only when the base field is locally compact. Notice that the only missing implication for a general non-Archimedean base field is (3) implies (2) or (3) implies (1).
\end{rem}

\begin{cor} \label{lem:compactoid_epi}
Let $f: E \to F$ be a morphism of nuclear Fr\'echet spaces. Then, the following are equivalent:
\begin{enumerate}
\item $f$ is an epimorphism in the category of complete locally convex vector spaces;
\item $f: E^b \to F^b$ is an epimorphism in $\ttCBorn_{k}$
\item $f: \ttCpt(E) \to \ttCpt(F)$ is an epimorphism in $\ttCBorn_{k}$;
\item $f: E^b \to F^b$ has bornologically dense image.
\end{enumerate}
\end{cor}

{\bf Proof.}
The characterization of epimorphisms given in Proposition \ref{prop:strict_morphisms_born_close} in combination with Lemma \ref{lem:compactoid_neumann} yields that conditions (2),  (3), and (4) are equivalent. Then, Lemma \ref{lem:compactoid_conv} implies that (1) is equivalent to (4) because Fr\'echet spaces are metric spaces. Therefore the topological closure of the image of $E$ in $F$ is equal to its limit points (for the topology of $F$). And since $F^b$ is a proper bornological vector space (because it is normal) we can apply Theorem \ref{prop:sequence_close} to deduce that the bornological closure of $E$ in $F$ agrees with its bornological limit point.
\ \hfill $\Box$

\begin{lem} \label{lem:short_exact_t}
If a short sequence of bornological Fr\'echet spaces
\[ 0 \to E \to F \to G \to 0 \]
is strictly exact in $\ttCBorn_{k}$ then 
\[ 0 \to E^t \to F^t \to G^t \to 0 \]
is strictly exact in $\ttLoc_k$.
\end{lem}

{\bf Proof.}
The ${(-)}^t$ functor is a left adjoint, so it preserve all colimits. Therefore, we have to check that it preserves kernels of morphisms of Fr\'echet spaces. Let's consider the kernel of the strict morphism
\[ f^t: F^t \to G^t. \]
$\ker(f^t)$ is a Fr\'echet space, and therefore it is a normal locally convex space. Hence, we can apply Proposition \ref{prop:limpro_commutation_t} to deduce that $\ker(f^t) \cong \ker(f)^t \cong E^t$, and the lemma is proved.
\ \hfill $\Box$

One problem that complicates our work is that the converse statement of Lemma \ref{lem:short_exact_t}  does not hold in general. However, we have the following result.

\begin{lem} \label{lem:compactoid_short_exact}
A short sequence of nuclear Fr\'echet spaces
\[ 0 \to E \to F \to G \to 0 \]
is strictly exact in $\ttLoc_k$ if and only if 
\[ 0 \to E^b \to F^b \to G^b \to 0 \]
is strictly exact in $\ttCBorn_k$.
\end{lem}

{\bf Proof.}
By Lemma \ref{lem:short_exact_t} the strict exactness of the second sequence implies the strict exactness of the first one. Applying the functor ${(-)}^b$  to the first sequence we obtain a strict exact sequence
\[ 0 \to E^b \to F^b \to G^b  \]
because ${(-)}^b$ is a right adjoint functor. Since the sequence
\[ 0 \to E^b \to F^b \to G^b \to 0 \]
is manifestly algebraically exact it remains to check that $F^b \to G^b$ is a strict epimorphism. By the nuclearity hypothesis on $E, F$ and $G$ we can apply Lemma \ref{lem:compactoid_neumann} to deduce that the von Neumann bornology of these spaces coincides with the compactoid one. Finally, we can apply Lemma \ref{lem:compactoid_quotient} to conclude the proof.
\ \hfill $\Box$

\begin{rem}
The functor ${(-)}^b$ does not preserve strict exactness of short exact sequences in general. It is known that there exist Montel-Fr\'echet spaces $V$ whose quotient $V/W$ by a closed subspace $W$ is not Montel. Since for such a $V$ one has $V^b \cong \ttCpt(V)$ whereas $(V/W)^b \not\cong \ttCpt(V/W)$, Lemma \ref{lem:compactoid_quotient} implies that 
\[ \ttCpt(V/W) \cong \frac{V^b}{W^b}. \]
So, ${(-)}^b$ does not preserve cokernels in this case.
\end{rem}

\subsection{Derived functors of the inverse limit functor} \label{sec:der_lim}

We assume the reader is familiar with the notion of a family of injective objects with respect to a functor between quasi-abelian categories, in the sense of Schneiders (Definition 1.3.2 of \cite{SchneidersQA}).
In this section we recall how to derive the inverse limit functors in quasi-abelian categories, and then focus on the case of $\ttCBorn_k$.  Conditions for the existence of the derived functor of the inverse limit functors are given the Proposition \ref{prop:inv_prosmans}, which is proven in Prosmans \cite{Pr3}.
We will discuss then a bornological version of the classical Mittag-Leffler Lemma for Fr\'echet spaces, which is an important example of the use of homological algebra towards functional analysis. There is an extensive literature on the study of the derived functor in the category of locally convex spaces: Palamadov \cite{Pal}, Vogt \cite{Vogt}, and Retakh \cite{Ret} are only some of the main contributions. In this section we bound ourself to obtain bornological versions of the basic results for Fr\'echet spaces, without looking for the utmost generality.

\begin{prop} \label{prop:inv_prosmans}
Let $I$ be a small category and $\ttC$ a quasi-abelian category with exact products. Then, the family of objects in $\ttC^{I^{op}}$ which are Roos-acyclic form a $\underset{i\in I}\limpro$-acyclic family. In particular, the functor 
\[\limpro_{i\in I}:\ttC^{I^{op}} \to \ttC \]
is right derivable to a functor 
\[D^{+}(\ttC^{I^{op}} ) \to D^{+}(\ttC )
\]
 and for any object $V \in \ttC^{I^{op}}$, we have a canonical isomorphism 
\begin{equation}\label{eqn:Derived2Roos}
\Rb \limpro_{i\in I} V_{i}\cong Roos(V)
\end{equation}
where the right hand side is the Roos complex of $V$.
\end{prop}
{\bf Proof.} 
See Section3.3 of \cite{Pr3}.
\hfill $\Box$

\begin{cor}The family of $\underset{i \in I}\limpro$-acyclic objects form a family of injective objects relative to the functor $\underset{i \in I}\limpro: \ttC^{I} \to \ttC$.
\end{cor}

In the case of a tower (\ie when $I = \Nb$ with its natural order), the situation is easier to deal with. Consider a functor $V:\mathbb{N}^{op}\to \ttC$ where $\ttC$ is a quasi-abelian category with exact products. We can consider the complex in degree $0$ and $1$

\begin{equation}\label{eqn:Weibel} \prod_{i=0}^{\infty} V_{i} \stackrel{\Delta}\longrightarrow \prod_{i=0}^{\infty} V_{i}
\end{equation}
defined by 
\[\Delta = (\dots, id_{V_2} - \pi_{2, 3}, id_{V_1} - \pi_{1, 2}, id_{V_0} - \pi_{0, 1})
\]
where $\pi_{i, j}: V_j \to V_i$ denote the canonical morphisms. In this particular case, the Roos complex reduces to the complex of equation (\ref{eqn:Weibel}).

\begin{lem}\label{lem:RoosML}
	Let $\ttC$ be a quasi-abelian category with exact products. Consider a functor $V:\mathbb{N}^{op}\to \ttC$. Then $\Rb \underset{i\in \mathbb{N}}\limpro V_{i}$ is isomorphic to the complex in degree $0$ and $1$ given in equation (\ref{eqn:Weibel}). 
\end{lem}
{\bf Proof.}
Since the cardinality of the natural numbers is less than the second infinite cardinal, Theorem 3.10 of \cite{Pr4} implies that $LH^{n}(\Rb \underset{i\in I}\limpro V_{i})=0$ for all $n\geq 2$. The fact that $\Rb \underset{i\in I}\limpro V_{i}$ is represented by (\ref{eqn:Weibel}) follows from the general definition of the Roos complex and the isomorphism in (\ref{eqn:Derived2Roos}). 
\hfill $\Box$

\begin{lem}\label{lem:TopML} Let $V$ be a projective system of Fr\'{e}chet spaces in $\ttLoc_k$ indexed by $\Nb$ where all system morphisms are dense. Then, the complex
\begin{equation} 
0 \to \limpro_{i \in \mathbb{N}}V_{i} \to \prod_{i \in \mathbb{N}}V_{i} \stackrel{\Delta}\to \prod_{i \in \mathbb{N}}V_{i} \to 0
\end{equation}
is strictly exact.
\end{lem}
{\bf Proof.} We can first show exactness in the category $\ttVect_{k}$ by applying the standard Mittag-Leffler theorem for (topological) Fr\'{e}chet spaces (see for example \cite{DV}, Lemma 1 of page 45) to 
\[\xymatrix{0 \ar[d]  & 0 \ar[l] \ar[d]& 0 \ar[l] \ar[d]& \cdots \ar[l] \ar[d] \\
V_1 \ar[d]  & V_2 \ar[l] \ar[d]& V_3 \ar[l] \ar[d]& \cdots \ar[l] \ar[d] \\
V_1 \ar[d] & V_1 \times V_2 \ar[l] \ar[d]& V_1 \times V_2 \times V_3\ar[d] \ar[l] & \cdots \ar[l] \ar[d]\\
0 \ar[d] & V_1 \ar[d] \ar[l]& V_1 \times V_2 \ar[d] \ar[l]& \cdots \ar[d] \ar[l] \\
0   & 0 \ar[l] & 0 \ar[l] & \cdots \ar[l] 
 .}
\] 
The lemma in \cite{DV} discusses only in the case $k = \Cb$, but it is easy to check the proof works over any base field since the only hypothesis used is that the spaces involved are endowed with a metric for which they are complete. It is an immediate consequence of the open mapping theorem for Fr\'echet spaces (see \cite{H2} Section4.4.7, page 61) that the morphism $\Delta$ is a strict epimorphism and therefore the sequence is strictly exact.
\hfill $\Box$

We can now state our bornological version of the Mittag-Leffler Lemma.

\begin{lem} (Mittag-Leffler) \\ \label{lemma:mittag_frechet}
Let $E, F, G \in \ttCBorn_k^{\Nb^{op}}$, be projective systems of bornological Fr\'echet spaces over $k$ indexed by $i \in \Nb$. Let 
\begin{equation} \label{eqn:short_ml}
0 \to \{ E_i \}_{i \in \Nb}  \stackrel{\eta}{\to}  \{ F_i \}_{i \in \Nb} \stackrel{\psi}{\to}  \{ G_i \}_{i \in \Nb} \to 0 
\end{equation} 
be a short exact sequence of systems where each $\eta_{i}$ and $\psi_{i}$ is strict in $\ttCBorn_k$. Suppose also, that 
\begin{enumerate}
\item $\{E_i \}_{i \in \Nb} $ is an epimorphic system;
\item $\underset{i \in \Nb}\limpro E_i$, $\underset{i \in \Nb}\limpro F_i$ and  $\underset{i \in \Nb}\limpro G_i$ are nuclear bornological spaces.
\end{enumerate}
Then, the resulting sequence 
\[ 0 \to \limpro_{i \in \Nb} E_i \to \limpro_{i \in \Nb} F_i \to \limpro_{i \in \Nb} G_i \to 0 \]
is strictly exact in $\ttCBorn_k$, where the limits are calculated in $\ttCBorn_k$.
\end{lem}

{\bf Proof.}
The datum of the short exact sequence (\ref{eqn:short_ml}) is equivalent to a sequence of strictly exact sequences
\[ i \mapsto (0 \to  E_i  \stackrel{\eta_{i}}\to F_i  \stackrel{\psi_{i}}\longrightarrow  G_i \to 0) \]
whose morphisms are compatible with the system morphisms of the projective systems. Then, the application of the functor ${(-)}^t$ to these sequences
\[ i \mapsto (0 \to  E_i^t  \stackrel{\eta_{i}}\to F_i^t  \stackrel{\psi_{i}}\longrightarrow  G_i^t \to 0) \]
yields strictly exact sequences in $\ttLoc_k$, thanks to Lemma \ref{lem:short_exact_t}. The system maps $E_{j}^{t} \to E_{i}^{t}$ have topologically dense set theoretic image by Observation \ref{obs:Dense2Dense}. Therefore, we can apply the Mittag-Leffler lemma for Fr\'echet spaces (cf. Lemma \ref{lem:TopML}) to the systems $\{E_i^t \}_{i \in \Nb} $ to get a strictly exact sequence in $\ttLoc_k$
\begin{equation} \label{eqn:short_ml2}
0 \to \limpro_{i \in \Nb} (E_i^t) \to \limpro_{i \in \Nb} (F_i^t) \to \limpro_{i \in \Nb} (G_i^t) \to 0
\end{equation} 
of Fr\'echet spaces. By Proposition \ref{prop:limpro_commutation_t}, the functors $\underset{i \in \Nb}\limpro$ and ${(-)}^t$ commute in (\ref{eqn:short_ml2}). Applying Lemma \ref{lem:compactoid_short_exact} we deduce that the strict exactness of the sequence 
\[ 0 \to (\limpro_{i \in \Nb} E_i)^t \to (\limpro_{i \in \Nb} F_i)^t \to (\limpro_{i \in \Nb} G_i)^t \to 0 \]
implies the strict exactness of the sequence
\[ 0 \to \limpro_{i \in \Nb} E_i \to \limpro_{i \in \Nb} F_i \to \limpro_{i \in \Nb} G_i \to 0, \]
concluding the proof.
\ \hfill $\Box$

\begin{cor}\label{cor:BornML} Let $V$ be a projective system of bornological Fr\'{e}chet spaces indexed by $\Nb$ where all system morphisms are dense and $\limpro V_i$ is nuclear. Then, the complex
	\begin{equation} 
	0 \to \limpro_{i \in \mathbb{N}}V_{i} \to \prod_{i \in \mathbb{N}}V_{i} \stackrel{\Delta}\to \prod_{i \in \mathbb{N}}V_{i} \to 0
	\end{equation}
	is strictly exact.
\end{cor}
{\bf Proof.}
Using the bornological version of the Mittag-Leffler Lemma, \ie Lemma \ref{lemma:mittag_frechet}, we can use the same proof of Lemma \ref{lem:TopML} to deduce the corollary.
\hfill $\Box$

\begin{cor}\label{cor:BornML2} Let $V$ be a projective system of bornological Fr\'{e}chet spaces indexed by $\Nb$ where all system morphisms are dense and $\limpro V_i$ is nuclear. Then
	\begin{equation} \label{eqn:ML_Roos}
	\Rb {\limpro_{i \in \mathbb{N}}} V_{i} \cong \limpro_{i \in \mathbb{N}}V_{i}.
	\end{equation}
\end{cor}
{\bf Proof.}
Corollary \ref{cor:BornML} is equivalent to say that the Roos complex of $V$ has cohomology only in degree $0$. Therefore, Proposition \ref{prop:inv_prosmans} implies (\ref{eqn:ML_Roos}).
\hfill $\Box$

\begin{lem}\label{lem:LongShort} 
In any quasi-abelian category $\ttC$, a complex 
\[\cdots \to V^{n+1} \stackrel{d^{n+1}}\to V^{n}  \stackrel{d^{n}}\to V^{n-1} \to \cdots \]
is strictly exact if and only if the sequences 
\begin{equation} \label{eqn:exact_sequence}
0 \to \ker(d^{n}) \to V^{n} \to \ker(d^{n-1}) \to 0
\end{equation}
are strictly exact for each $n$. Given a projective system of strictly exact complexes 
\begin{equation} \label{eqn:exact_complex}
\cdots \to V_{i}^{n+1} \stackrel{d_{i}^{n+1}}\to V_{i}^{n}  \stackrel{d_{i}^{n}}\to V_{i}^{n-1} \to \cdots
\end{equation}
where both  $\{\ker(d^{n}_{i})\}_{i \in \Nb}$ and $\{V_{i}^{n}\}_{i \in I}$ are $\underset{i \in \Nb}\limpro$-acyclic systems, the projective limit
\[\cdots \to \limpro_{i \in \Nb}V_{i}^{n+1} \stackrel{\limpro d_{i}^{n+1}}\to \limpro_{i \in \Nb}V_{i}^{n}  \stackrel{\limpro d_{i}^{n}}\to\limpro_{i \in \Nb} V_{i}^{n-1} \to \cdots \] 
is strictly exact.
\end{lem}
{\bf Proof.}
By Corollary 1.2.20 of \cite{SchneidersQA}, the complex (\ref{eqn:exact_complex}) is strictly exact if and only if $LH^{n}(V)=0$ for all $n$. This is equivalent to the two term complexes $0 \to \coim(d^{n + 1}) \to \ker(d^n) \to 0$ vanishing in the left heart of $\ttC$ (which is a full subcategory of the derived category, cf. Corollary 1.2.20 of ibid.). By this explicit description of the left heart, the condition $LH^{n}(V)=0$ is equivalent to that the canonical morphism $\coim(d^{n + 1}) \to \ker(d^n)$ is an isomorphism. Therefore, the strictly exactness of sequences
\[ 0 \to \ker(d^{n}) \to V^{n} \to \coim(d^{n}) \to 0 \]
implies the strictly exactness of the sequences (\ref{eqn:exact_sequence}). For the second statement, observe that the sequences 
\[0 \to \ker{d^{n}} \to  \limpro_{i\in \mathbb{N}} V^{n}_{i} \to \ker{d^{n-1}} \to 0 \]
are strictly exact, being the application of $\mathbb{R} \underset{i\in \mathbb{N}}\limpro$ to the strict short exact sequences 
\[0 \to \ker{d_{i}^{n}} \to  V^{n}_{i} \to \ker{d_{i}^{n-1}} \to 0 \]
of $\underset{i\in \mathbb{N}}\limpro$-acyclic objects, thought of as an exact triangle in $D^{+}(\ttC^{\mathbb{N}})$.
\hfill $\Box$

\section{Stein domains} \label{sec:Stein}

\subsection{Bornological Fr\'echet algebras}

The notion of multiplicatively convex bornological algebra (in short m-algebra) introduced in \cite{BaBe} (cf. Definition 4.1 of ibid.) is not general enough for all purposes of analytic geometry. For example, the bornological Fr\'echet algebras of analytic functions on open subsets of $\Ab_k^n$ are not multiplicatively convex in this sense. So, we introduce here a generalization of multiplicatively convex bornological algebras which encompass also bornological Fr\'echet algebras. We start by recalling what the spectrum of a general bornological algebra is.

\begin{defn}
	Let $A$ be a bornological algebra, \ie an object of $\ttComm(\ttBorn_k)$, we define the \emph{spectrum} of $A$ as the set of equivalence classes of bounded algebra morphisms of $A$ to valued extensions of $k$. The spectrum is denoted by $\mM(A)$ and it is equipped with the weak topology, \ie the weakest topology for which all maps of the form $\chi \mapsto |\chi(f)| \in \Rb_{\ge 0}$, for all $f \in A$, are continuous.
\end{defn}

This definition extends the one given in Definition 4.2 of \cite{BaBe} for m-algebras, but the spectrum of a general bornological algebra is not as well behaved as the spectrum of a bornological m-algebra. For example, there exist bornological algebras whose underlying bornological vector space is complete and whose spectrum is empty. So, it is important to single out a suitable sub-category of the category of bornological algebras for which the notion of spectrum is not pathological.

\begin{defn}\label{defn:ProComplete}
	Let $\A$ be an object of $\ttComm(\ttCBorn_{k})$. We say that $\A$ is \emph{pro-multiplicatively convex} (or a \emph{pro m-algebra}) if there is an isomorphism
	\[ \A \cong \limpro_{i \in I} \A_i \]
	in $\ttComm(\ttBorn_{k})$, where $I$ is a cofiltered small category and $\A_i$ are complete bornological m-algebras.
\end{defn}

\begin{defn}\label{defn:DenseStrictlyDense}
	Let $\A $ be a pro-multiplicatively convex object of $\ttComm(\ttCBorn_{k})$. $A$ is called \emph{densely defined} if there exists an isomorphism $\A \cong \underset{i \in I}\limpro \A_i$ of bornological algebras as in Definition \ref{defn:ProComplete} such that for any $i \in I$ the canonical map $\pi_i:\A \to \A_i$ has dense set-theoretic image.
\end{defn}

We remark that the condition in Definition \ref{defn:DenseStrictlyDense} can be stated purely categorically by requiring that the morphisms $\pi_i: \A \to \A_i$ are epimorphisms of the underlying complete bornological vector spaces (cf. Proposition \ref{prop:strict_morphisms_born_close}).

\begin{prop} \label{prop:frechet_spectrum}
	Let $A$ be a densely defined pro-multiplicatively convex bornological algebra, such that
	\[ A \cong \limpro_{n \in \Nb} A_n \]
	with $A_n$ Banach algebras. Then, $A$ is a bornological Fr\'echet algebra whose spectrum coincides with the spectrum of $\mM(A^t)$, \ie
	\[ \mM(A) = \bigcup_{n \in \Nb} \mM(A_n) \]
	topologically.
\end{prop}

{\bf Proof.}
$A$ is a bornological Fr\'echet algebra, \ie the underlying bornological vector space of $A$ is that of a bornological Fr\'echet space. Indeed, $\underset{n \in \Nb}\limpro (A_n^t)$ is a Fr\'echet space as consequence of Proposition \ref{prop:limpro_commutation_t}, which states that $\underset{n \in \Nb}\limpro (A_n^t) \cong (\underset{n \in \Nb}\limpro A_n)^t = A^t$. Thus,
\[ (A^t)^b \cong \limpro_{n \in \Nb} ((A_n^t)^b) \cong A. \]
This implies that a character $A \to K$ is bounded if and only if it is continuous for $A^t$. The  densely defined hypothesis implies that 
\[ \mM(A) \cong \mM(A^t) = \bigcup_{n \in \Nb} \mM(A_n) \]
because $\pi_i: A \to A_i$ is an epimorphism if and only if $\pi_i^t: A^t \to A_i^t$ is an epimorphism, since the $A_i$ are Banach. The identification
\[ \mM(A^t) = \bigcup_{n \in \Nb} \mM(A_n) \]
is a well-known result of the classical theory of Frech\'et algebras (\eg, see Section2.5 of \cite{Bam}).
\ \hfill $\Box$

We end this section recalling a result from \cite{Bam2} that will be used in the next subsection.

\begin{defn} \label{defn:born_web}
Let $E$ be a separated bornological vector space of convex type over $k$. A pair $(\mV, b)$ consisting of mappings $\mV : \underset{j \in \Nb}\bigcup \Nb^j \to \mathcal{P}(E)$ and $b : \Nb^\Nb \to (|k^\times|)^\Nb$ is called a \emph{bornological web} if all of the conditions below hold:
\begin{enumerate}
\item The image of $\mV$ consists of disks.
\item $\mV(\void) = E$.
\item Given a finite sequence $(n_0, \dots, n_j)$, then $\mV(n_0, \dots, n_j)$ is absorbed by
      \[ \bigcup_{n \in \Nb} \mV(n_0, \dots, n_j, n). \]
\item For every $s: \Nb \to \Nb$ the series $\underset{n \in \Nb}\sum \la(s)_n x_n$, with $\la(s)_n \in k$, converges bornologically in $E$, whenever we choose $x_n \in \mV(s(0), \dots ,s(n))$ and $|\la(s)_n| = b(s)_n$.
\end{enumerate}
\end{defn}

\begin{lem} \label{lem:closed_graph}
	Let $A, B$ be objects of $\ttComm(\ttBorn_{k})$ for which the underlying bornological vector space of $A$ is complete and the one of $B$ is a webbed bornological vector space. Let $\phi: A \to B$ be a morphism of the underlying objects in $\ttComm(\ttVect_{k})$. Suppose that in 
	$B$ there is a family of ideals $\Im$ such that
	\begin{enumerate}
		\item each $I \in \Im$ is (bornologically) closed in $B$ and each $\phi^{-1}(I)$ is closed in $A$;
		\item for each $I \in \Im$ one has $\dim_k B/I < \infty$;
		\item $\underset{I \in \Im} \bigcap I = (0)$.
	\end{enumerate}
	Then, $\phi$ is bounded.
\end{lem}

{\bf Proof.}
The proof can be found in \cite{Bam2}, Proposition 4.23.
\ \hfill $\Box$

And finally the following lemma permits us to apply Lemma \ref{lem:closed_graph} to the morphisms of bornological algebras we will study next sub-section.

\begin{lem} \label{lem:webbed}
	The underlying bornological vector space of every $k$-dagger affinoid algebra is a webbed bornological vector space.
\end{lem}

{\bf Proof.}
The assertion about dagger affinoid algebras is a direct consequence of Example 2.3 (2) of \cite{Bam2}.

\ \hfill $\Box$

\begin{rem}
	The notion of pro-multiplicative algebra is more general than what strictly needed in this paper. However, this notion is very comfortable when both Fr\'echet algebras and LB algebras are considered in the same discussion. Moreover, the material in this section will be a reference for future works in which we will analyse dagger quasi-Stein algebras for which the use of the notion of pro-multiplicative bornological algebra is unavoidable.
\end{rem}

\subsection{Stein algebras and Stein spaces} 

\begin{defn}
A \emph{Weierstrass localization} of a $k$-dagger affinoid algebra $\A$ is a morphism of the form
\[ \A \to \frac{ \A \lt r_1^{-1} X_1, \ldots, r_n^{-1} X_n \gt^\dagger}{(X_1 - f_1, \ldots, X_n - f_n)} \]
for some $f_1, \dots, f_n \in \A$, $(r_i) \in \Rb_+^n$. 
\end{defn}

\begin{defn}
	Let $A, B, C$ be a $k$-dagger affinoid algebras such that $B$ and $C$ are $A$-algebras. Let $f: B \to C$ be a bounded morphism of $A$-algebras. $f$ is called \emph{inner with 
		respect to $A$} if there exists a strict epimorphism $\pi: A \lt r_1^{-1} T_1, \dots, r_n^{-1} T_n \gt^\dagger \to B$ such that
	\[ \rho_C(f(\pi(T_i))) < r_i \]
	for all $1 \le i \le n$, where $\rho_{C}$ is the spectral semi-norm of $C$.
\end{defn}

\begin{defn}
	Let $\phi: \mM(A) = X \to \mM(B) = Y$ be a morphism of $k$-dagger affinoid spaces. The \emph{relative interior of $\phi$} is the set 
	\[ \ttInt(X/Y) = \{ x \in X | A \to \mH(x) \text{ is inner with respect to $B$} \}. \]
	The complement of $\ttInt(X/Y)$ is called the \emph{relative boundary of $\phi$} and denoted by $\partial(X/Y)$. If $B = k$, the sets $\ttInt(X/Y)$ and $\partial(X/Y)$ are denoted by $\ttInt(X)$ and $\partial(Y)$ and called the \emph{interior} and the \emph{boundary} of $X$.
\end{defn}

\begin{defn} \label{defn:stein_algebra}
	A \emph{dagger Stein algebra} over $k$ is a complete bornological algebra $A$ over $k$ which is isomorphic to an inverse limit of $k$-dagger affinoid algebras 
	\begin{equation} \label{eqn:stein}
	\cdots \longrightarrow A_4 \longrightarrow A_3 \longrightarrow A_2 \longrightarrow A_1 \longrightarrow A_0
	\end{equation} 
	in the category $\ttCBorn_{k}$ where each morphism is a Weierstrass localization and $\mM(A_i)$ is contained in the interior of $\mM(A_{i+1})$, for each $i$. The category of dagger Stein algebras is the full sub-category of $\ttComm(\ttCBorn_{k})$ identified by dagger Stein algebras over $k$ and it is denoted by $\ttStn_k$.
\end{defn}

\begin{defn} \label{defn:stein_spaces}
	A $k$-dagger analytic space $X$ is called \emph{dagger Stein space} if it admits a dagger affinoid covering $U_1 \subset U_2 \subset \cdots $ such that
	\[ X = \bigcup_{i \in \Nb} U_i \]
	and the restriction morphisms $\mO_X(U_{i+1}) \to \mO_X(U_i)$ are Weierstrass localizations and $U_i$ is contained in the interior of $U_{i+1}$, for each $i \in \Nb$. The category of dagger Stein spaces is the full sub-category of the category of $k$-analytic spaces of this form.
\end{defn}

\begin{example} \label{exa:stein}
\begin{enumerate}
\item The open polycylinders are the most basic example of Stein spaces, whose exhaustion by Weierstrass subdomains is given by the closed polycylinders of smaller radius. 
\item Analytifications of algebraic varieties are Stein spaces, both for Archimedean and non-Archimedean base fields.
\item Closed subspaces of Stein spaces are Stein spaces (see \cite{KI}).
\end{enumerate}
\end{example}

\begin{lem} \label{lemma:regular_quasi_Stein}
	Let $A$ be a dagger Stein algebra over $k$. Then, $A$ is a densely defined, pro-multiplicatively convex, bornological Fr\'echet algebra over $k$.
\end{lem}

{\bf Proof.}
Let $\underset{i \in \Nb}\limpro A_i$ be a presentation of $A$ as in definition of dagger Stein algebra. Let $\tilde{A}_i$ be the completion of $A_i$ with respect to the norm 
\[ \|f\| = \inf_{B \in \mD_{A_i}} |f|_{(A_i)_B}, \]
where $|\cdot|_{(A_i)_B}$ denotes the gauge semi-norm, as explained in Remark \ref{rem:gauge}.
Notice that $\|\cdot\|$ is a norm because by hypothesis $U_i$ has non-empty interior, therefore $\|f\| = 0 \then f = 0$ (more precisely, since $U_i$ has non-empty interior one can find an injective bounded morphism of $A_i$ to a Banach algebra, which implies that $A_i$ is separated). The requirement of $\mM(A_i)$ to lie in the interior of $\mM(A_{i+1})$ can be restated by saying that there are (bounded) diagonal morphisms such that the following diagram is commutative
\begin{equation} \label{eqn:diag_A_A_tilde} 
\begin{tikzpicture}
\matrix(m)[matrix of math nodes,
row sep=2.6em, column sep=2.8em,
text height=1.5ex, text depth=0.25ex]
{ \cdots  &  A_3 & A_2 & A_1  \\
	\cdots  &  \tilde{A}_3 & \tilde{A}_2 & \tilde{A}_1 &    \\};
\path[->,font=\scriptsize]
(m-1-1) edge node[auto] {$$} (m-1-2);
\path[->,font=\scriptsize]
(m-1-2) edge node[auto] {$$} (m-1-3);
\path[->,font=\scriptsize]
(m-1-3) edge node[auto] {$$} (m-1-4);
\path[->,font=\scriptsize]
(m-2-1) edge node[auto] {$$} (m-2-2);
\path[->,font=\scriptsize]
(m-2-2) edge node[auto] {$$} (m-2-3);
\path[->,font=\scriptsize]
(m-2-3) edge node[auto] {$$} (m-2-4);
\path[->,font=\scriptsize]
(m-1-2) edge node[auto] {$$} (m-2-2);
\path[->,font=\scriptsize]
(m-1-3) edge node[auto] {$$} (m-2-3);
\path[->,font=\scriptsize]
(m-1-4) edge node[auto] {$$} (m-2-4);
\path[->,font=\scriptsize]
(m-2-2) edge node[auto] {$$} (m-1-3);
\path[->,font=\scriptsize]
(m-2-3) edge node[auto] {$$} (m-1-4);
\end{tikzpicture}
\end{equation} 
where the vertical maps are the completions. Hence, following the diagonal maps we obtain a projective system
\begin{equation} \label{eqn:system_A_A_tilde} 
\cdots \to A_3 \to \tilde{A}_3 \to A_2 \to \tilde{A}_2 \to A_1 \to \tilde{A}_1, 
\end{equation} 
which shows that the sub-system
\[  \cdots \to \tilde{A}_3 \to \tilde{A}_2 \to \tilde{A}_1 \]
is cofinal in (\ref{eqn:system_A_A_tilde}), and hence
\begin{equation} \label{eqn:frechet_arens}
A \cong \limpro_{i \in \Nb} \tilde{A}_i. 
\end{equation} 
Then, by Proposition \ref{prop:limpro_commutation_t} we have $\underset{i \in \Nb}\limpro (\tilde{A}_i^t)  \cong (\underset{i \in \Nb}\limpro (\tilde{A}_i))^{t}$, because a countable projective limit of Banach spaces in $\ttLoc_k$ is a Fr\'echet spaces and hence normal (see Example \ref{example:normal}). So, from the relation 
\[A= \underset{i \in \Nb}\limpro (\tilde{A}_i) \cong (\underset{i \in \Nb}\limpro (\tilde{A}_i))^{tb} \cong (\underset{i \in \Nb}\limpro (\tilde{A}_i^t))^{b}
\]
we can deduce that $A$ is a Fr\'echet space in $\ttBorn_k$, in the sense of Definition \ref{defn:frechet}. Moreover, it is easy to see all the maps of the system (\ref{eqn:frechet_arens}) have dense images because in the diagram (\ref{eqn:diag_A_A_tilde}) all maps but the bottom horizontal ones are clearly epimorphisms, which implies that the bottom horizontal ones are epimorphism too.
\hfill $\Box$

\begin{rem} \label{rem:strict_affinoid}
For non-Archimedean base fields, we could have defined Stein algebras using classical affinoid algebras in place of dagger affinoid algebras (which is what we get in Lemma \ref{lemma:regular_quasi_Stein} with the algebras $\tilde{A}_i$). The cofinality argument of Lemma \ref{lemma:regular_quasi_Stein} shows that the bornological structure we defined on a dagger Stein algebra agrees with the Fr\'echet structures defined classically, for example in \cite{Ber1990}, \cite{Ch}, \cite{Cr} for the non-Archmedean theory and in \cite{GR} for the Archimedean theory. 
\end{rem}

\begin{lem} \label{lemma:stein_nuclear}
Let $A$ be a dagger Stein algebra over $k$. Then, the underlying bornological vector space of $A$ is binuclear.
\end{lem}

{\bf Proof.}
By Proposition \ref{prop:dagger_nuclear} the underlying bornological vector spaces of dagger affinoid algebras are nuclear. We can apply Proposition \ref{prop:nuclear_proj} to any system of the form (\ref{eqn:stein}) to deduce that the underlying bornological vector space of $A$ is nuclear. It remains to see that the locally convex space $A^t$ is nuclear. Notice that if we chose a presentation $A \cong \underset{i \in \Nb}\limpro A_i$, we have that $A^t \cong \underset{i \in \Nb}\limpro (A_i^t)$, as a consequence of Proposition \ref{prop:limpro_commutation_t}, and by Lemma \ref{lem:LB_nuclear} $A_i^t$ are nuclear locally convex spaces. This implies that $A^t$ is nuclear because projective limits of nuclear locally convex spaces are always nuclear (for a proof of this fact we can refer to \cite{DV} Proposition 3.2 at page 23 in the case in which $k$ is Archimedean, and for the non-Archimedean case to Theorem 8.5.7 of \cite{PGS}).
\hfill $\Box$

\begin{lem} \label{lemma:locally_compact_spectrum}
Let $A$ be a dagger Stein algebra over $k$. Then, $\mM(A)$ is a dagger Stein space over $k$.
\end{lem}

{\bf Proof.}
The spectrum $\mM(A)$ is a hemi-compact topological space because, by Lemma \ref{lemma:regular_quasi_Stein}, we can apply Proposition \ref{prop:frechet_spectrum} to deduce that $\mM(A) = \underset{i \in \Nb}\bigcup \mM(A_i)$ (topologically), for a presentation $A \cong \underset{i \in \Nb}\limpro A_i$ (because, in the notation of Lemma \ref{lemma:regular_quasi_Stein}, $\mM(A_i) \cong \mM(\tilde{A}_i)$). The condition of $\mM(A_i)$ being contained in the Berkovich interior of $\mM(A_{i+1})$ readily implies that $\{\mM(A_i) \}_{i \in \Nb}$ is a cofinal family of compact subsets of $\mM(A)$ (cf. Exercise 3.8.C (b) of \cite{Eng}). This family is therefore a Berkovich net on $\mM(A)$ which induces a structure of $k$-dagger analytic space. It is straightforward to check that in this way $\mM(A)$ is endowed with a structure of dagger Stein space because by the definition of dagger Stein algebra the morphisms $A_{i + 1} \to A_i$ are Weierstrass localizations and $\mM(A_i)$ are mapped into $\ttInt(\mM(A_{i + 1}))$ by these localizations.
\hfill $\Box$

\begin{lem} \label{lemma:Stein_algebra_strict}
Let $A$ be a dagger Stein algebra over $k$. Then, we can always find a presentation
\[ A \cong \limpro_{i \in \Nb} A_i \]
as in Definition \ref{defn:stein_algebra} such that $A_i$ are strictly dagger affinoid algebras.
\end{lem}
{\bf Proof.}
Let $A \cong \underset{i \in \Nb}\limpro A_i$ be a presentation of $A$ as in the definition of dagger Stein algebra.
It is enough to show that each localization $A_{i + 1} \to A_i$ factors through Weierstrass localizations $A_{i + 1} \to B_i \to A_i$ where $B_i$ are strictly dagger affinoid. It is easy to check that Lemma 2.5.11 of \cite{Ber1990} holds also for dagger affinoid algebras. Therefore, for any $0 < \epsilon < 1$ we can find a $\rho = (\rho_i) \in \Rb_+^n$ and a strict epimorphism
\[ \pi: k \lt \rho_1^{-1} X_1, \ldots, \rho_n^{-1} X_n \gt^\dagger \to A_{i + 1} \]
such that the image of $\mM(A_i)$ is contained in the image of $\mM(A \lt (\epsilon \rho_1)^{-1} \pi(X_1), \ldots, (\epsilon \rho_n)^{-1} \pi(X_n) \gt^\dagger)$ in $\mM(A_{i + 1})$. Since the diagram
\[
\begin{tikzpicture}
\matrix(m)[matrix of math nodes,
row sep=2.6em, column sep=2.8em,
text height=1.5ex, text depth=0.25ex]
{   k \lt \rho_1^{-1} X_1, \ldots, \rho_n^{-1} X_n \gt^\dagger & A_{i + 1}  \\
	k \lt (\epsilon \rho_1)^{-1} X_1, \ldots, (\epsilon \rho_n)^{-1} X_n \gt^\dagger & A_{i + 1} \lt (\epsilon \rho_1)^{-1} \pi(X_1), \ldots, (\epsilon \rho_n)^{-1} \pi(X_n) \gt^\dagger \\};
\path[->,font=\scriptsize]
(m-1-1) edge node[auto] {$\pi$} (m-1-2);
\path[->,font=\scriptsize]
(m-2-1) edge node[auto] {$$} (m-2-2);
\path[->,font=\scriptsize]
(m-1-1) edge node[auto] {$$} (m-2-1);
\path[->,font=\scriptsize]
(m-1-2) edge node[auto] {$\pi'$} (m-2-2);
\end{tikzpicture}
\]
is a push-out square, we see that the morphism $\pi'$ is a strict epimorphism, because tensoring preserves strict epimorphisms. Therefore, choosing $(\epsilon \rho_i)^{-1} \in \sqrt{|k^\times|}$ for all $i$, which is always a possible choice because $\sqrt{|k^\times|}$ is dense in $\Rb$, we see that $A_{i + 1} \lt (\epsilon \rho_1)^{-1} \pi(X_1), \ldots, (\epsilon \rho_n)^{-1} \pi(X_n) \gt^\dagger$ can always be found to be a strictly dagger affinoid algebra and this is our choice for $B_i$. Moreover, the morphism $A_{i + 1} \to B_i$ is by construction a Weierstrass subdomain embedding and $A_{i + 1} \to A_i$ is by hypothesis a Weierstrass subdomain embedding which immediately implies that $B_i \to A_i$ is a Weierstrass domain embedding, concluding the proof.
\hfill $\Box$

\begin{rem}
In Lemma \ref{lemma:Stein_algebra_strict} we only discussed how to find a presentation of a dagger Stein algebra as an inverse limit of strictly dagger affinoid algebras, but the same argument clearly works in the classical affinoid case. Notice that this lemma crucially depends on the hypothesis that $k$ is non-trivially valued.
\end{rem}

\begin{lem} \label{lemma:quasi_Stein_algebra_closed_ideals}
Let $A$ be a dagger Stein algebra over $k$ and $\fm \subset A$ a finitely generated maximal ideal. Then, $\fm$ is bornologically closed in $A$.
\end{lem}
{\bf Proof.}
Let $A \cong \underset{i \in \mathbb{N}}\limpro A_i$ be a presentation of $A$ as given by Lemma \ref{lemma:Stein_algebra_strict}. Consider the extensions of $\fm \subset A$ to $A_i$, denoted $\fm A_i$, for every $i \in \Nb$. There are two possible cases to consider: either $\fm A_i$ is a proper ideal of $A_i$, for some $i$, or $\fm A_i = A_i$ for all $i$.  In the first case, since $\fm$ is maximal in $A$ then $\pi_i^{-1}(\fm A_i)$ must be a proper ideal containing $\fm$, and therefore must be equal to $\fm$. Since $\fm A_i$ is bornologically closed in $A_i$ (cf. Theorem 4.9 (2) of \cite{BaBe}), $\pi_i$ is a bounded map and the pre-image of a bornologically closed set by bounded maps is a bornologically closed set, then $\fm$ is bornologically closed in $A$.

On the other hand, suppose that for every $i$, $\fm A_i = A_i$ and let $f_1, \ldots, f_n$ denote a set of generators of $\fm$. Since $\fm A_i = A_i$, $\pi_i(f_i)$ have no common zeros in $\mM(A_i)$. Because this is true for any $i$ and $\mM(A) = \underset{i \in \mathbb{N}}\bigcup \mM(A_i)$, we can deduce that $f_1, \ldots, f_n$ have no common zeros in $\mM(A)$. Using Theorem A for dagger Stein spaces (cf. Theorem 3.2 of \cite{GK} for the case in which $k$ is non-Archimedean, instead for $k$ Archimedean the classical theorem of Cartan applies) and reasoning in the same way as Theorem V.5.4 and Theorem V.5.5 at page 161 of \cite{GR}, we can deduce that condition that $f_1, \ldots, f_n$ have no common zeros implies that there exist $g_1, \ldots, g_n \in A$ such that $1 = \underset{i = 1}{\overset{n}\sum} f_i g_i$ and hence $\fm = A$. But this is impossible because by hypothesis $\fm$ is a proper ideal of $A$. This contradiction shows that there must exist an $i$ such that $\fm A_i \ne A_i$.
\hfill $\Box$

\begin{rem}
The second part of the proof of Lemma \ref{lemma:quasi_Stein_algebra_closed_ideals} works only for fintely generated ideals. Indeed, dagger Stein algebras can have non-finitely generated maximal ideals. It follows that such ideals are necessarily non-closed and bornologically dense subsets.
\end{rem}

\begin{rem} \label{rem:immersion_A_i_A}
Notice that the morphism of Grothendieck locally ringed spaces $\mM(A_i) \to \mM(A)$, induced from the projections $\pi_i: A \to A_i$ discussed so far, is an open immersion (in the sense of Definition 4.18 of \cite{Bam}), when one consider the Berkovich G-topology of analytic domains (cf. chapter 6 of \cite{Bam} or chapter 1 of \cite{Ber1993} for details about this G-topology).
This is because we endow $\mM(A)$ with the structure of $k$-dagger analytic space given by the Berkovich net $\underset{i \in \Nb}\bigcup \mM(A_i)$. Therefore, each $\mM(A_i)$ is an analytic domain in $\mM(A)$ and $\mM(A_i)$ is identified with its image in $\mM(A)$.
\end{rem}

\begin{lem}\label{lem:FactorizationQSt}
Any morphism of $k$-dagger Stein algebras comes from a morphism in $\ttPro(\ttAfnd^{\dagger}_{k})$ by applying the functor \[\limpro :\ttPro(\ttAfnd^{\dagger}_{k}) \to \ttComm(\ttBorn_{k}).\]
\end{lem}
{\bf Proof.}
Let $f: A \to B$ be a morphism of dagger Stein algebras and let $A \cong \underset{i \in \mathbb{N}}\limpro A_i$, $B \cong \underset{j \in \mathbb{N}}\limpro B_j$ be two fixed presentations as in the definition of dagger Stein algebra. Let $\mM(f): \mM(B) \to \mM(A)$ denote the morphism of dagger Stein spaces induced by $f$. 
Since the family $\{ \mM(A_i)\}_{i \in \Nb}$ is cofinal in the family of compact subsets of $\mM(A)$, every morphisms $\mM(B_i) \to \mM(A)$ must factor through a $\mM(A_i)$, because the image of $\mM(B_i)$ in $\mM(A)$ is compact.  Therefore, we get morphisms of $k$-dagger analytic spaces $\mM(f_i): \mM(B_i) \to \mM(A_i)$ (possibly by re-indexing the systems in a suitable way). These morphisms commute with the morphisms $\mM(f): \mM(B) \to \mM(A)$ and $\underset{i \in \Nb}\limind \mM(f_i) = \mM(f)$. By Remark 4.16 of \cite{BaBe}, there exist maps $f_i: A_i \to B_i$ of $k$-dagger affinoid algebras induced by the morphisms of $k$-dagger affinoid spaces $\mM(f_i)$ and for these maps we get $\underset{i \in \Nb}\limpro f_i = f$, as required. 
\ \hfill $\Box$

\begin{thm} \label{prop:quasi_Stein_algebra_spaces}
The category of dagger Stein algebras over $k$ is anti-equivalent to the category of dagger Stein spaces over $k$. Moreover, the forgetful functor from the category of dagger Stein algebras over $k$ to $\ttComm(\ttVect_{k})$ is fully faithful.
\end{thm}
{\bf Proof.}
It is clear that to any dagger Stein algebra one can associate a dagger Stein space and vice versa functorially, as explained so far (cf. Lemma \ref{lemma:locally_compact_spectrum}). So, we check the claim about the boundedness of every algebra morphism between dagger Stein algebras. Fix two dagger Stein algebras with representations $A \cong \underset{i \in \Nb} \limpro A_i$, $B \cong \underset{i \in \Nb}\limpro B_i$, and let $\phi: A \to B$ be a morphism of their underlying objects in $\ttComm(\ttVect_{k})$. We also suppose that all $A_i$ and $B_i$ are strictly dagger affinoid algebras, which is always a possible choice (cf. Lemma \ref{lemma:Stein_algebra_strict}).
The data of the morphism $\phi$ is equivalent to a system of morphisms $\phi_i: A \to B_i$ in $\ttComm(\ttVect_{k})$ with obvious commuting relations with the system morphism of the representation of $B$ we fixed. To show that $\phi$ is bounded, it is enough to show that each $\phi_i$ is bounded. Let $\fm \subset B_i$ be a maximal ideal. Then, $\fm$ is finitely generated and $\frac{B_i}{\fm}$ is a finite valued extension of $k$, because $B_i$ is strictly dagger affinoid. Then, the composition of maps
\[ A \to B_i \to \frac{B_i}{\fm} \]
identifies $\frac{A}{\phi_i^{-1}(\fm)}$ with a sub-ring of $\frac{B_i}{\fm}$ so 
\[k \subset \frac{A}{\phi_i^{-1}(\fm)} \subset \frac{B_i}{\fm}.\] 
This implies of course that $\frac{A}{\phi_i^{-1}(\fm)}$ is a field. Therefore $\phi_i^{-1}(\fm)$ is a maximal, finitely generated ideal of $A$, which by Lemma \ref{lemma:quasi_Stein_algebra_closed_ideals} must be closed. It follows that the quotient bornology of $\frac{A}{\phi_i^{-1}(\fm)}$ is separated, and hence complete and since $\frac{A}{\phi_i^{-1}(\fm)}$ as vector space over $k$ is finite dimensional, its bornology is necessarily isomorphic to the product bornology of finitely many copies of $k$. Lemma \ref{lemma:locally_compact_spectrum} yields
\[ \mM(A) = \bigcup_{i\in \mathbb{N}} \mM(A_i) \]
which readily implies
\[ \Max(A) = \bigcup_{i\in \mathbb{N}} \Max(A_i), \]
where on the left-hand side only finitely generated maximal ideals are considered. So, there exists a $j$ such that for any $j \ge i$, $\phi_i^{-1}(\fm)A_j$ is a maximal ideal of $A_j$. It follows that there is a canonical isomorphism
\[ \frac{A}{\phi_i^{-1}(\fm)} \cong \frac{A_j}{\phi_i^{-1}(\fm)A_j}. \]
Now consider $\fm^n$, for any $n \in \Nb$. By elementary commutative algebra one can see that $\phi_i(\fm)^n \subset \phi_i(\fm^n)$ which implies that there exist canonical quotient maps
\[ A \to \frac{A}{\phi_i^{-1}(\fm)^n} \to \frac{A}{\phi_i^{-1}(\fm^n)}. \]
Since $\mM(A_j) \rhook \mM(A)$ is an immersion (\ie induces isomorphisms on stalks), using the same argument of proposition 7.2.2/1 of \cite{BGR} we can deduce that
\[ \frac{A}{\phi_i^{-1}(\fm)^n} \cong \frac{A_j}{\phi_i^{-1}(\fm)^n A_j}, \ \ \forall n \in \Nb. \]
Notice also that $\frac{A_j}{\phi_i^{-1}(\fm)^n A_j}$ are a finite dimensional $k$-Banach algebra. Therefore, we wrote $\frac{A}{\phi_i^{-1}(\fm^n)}$ is a quotient of a finite dimensional $k$-Banach algebra, which is necessarily a $k$-Banach. This shows that $\phi_i^{-1}(\fm^n)$ is closed in $A$, for any maximal ideal of $B_i$ and any $n \in \Nb$. We showed that the family of all powers of maximal ideals of $B_i$ satisfies the hypothesis of Lemma \ref{lem:closed_graph} (that can be applied to $\phi_i$ thanks to Lemma \ref{lem:webbed}), proving that $\phi_i$ and hence $\phi$ are bounded maps.

To conclude the proof it remains to show that every morphism of dagger Stein spaces is induced by a morphism of dagger Stein algebras. A morphism  $X \to Y$ of dagger Stein spaces is by definition a morphism of locally ringed Grothendieck topological spaces between $X = \underset{i \in \mathbb{N}}\bigcup \mM(B_i)$ and $Y= \underset{i\in \mathbb{N}}\bigcup \mM(A_i)$, when they are equipped with their maximal Berkovich nets generated by the nets $\{ \mM(A_i) \}_{i \in \Nb}$ and $\{\mM(B_i)\}_{i \in \Nb}$ (cf. chapter 6 of \cite{Bam}). All the elements of the maximal nets which gives the structure of $k$-analytic spaces to $X$ and $Y$ must be compact subsets of $X$ and $Y$ respectively, and hence by the fact that $X$ and $Y$ are hemi-compact, all the elements of the maximal nets must be contained in some $\mM(A_i)$ or $\mM(B_i)$, respectively, as explained in the proof of Lemma \ref{lemma:locally_compact_spectrum}. Therefore, every morphism of $k$-dagger analytic spaces must come from a morphism of systems $i \mapsto (A_{i} \to B_{i})$, of the form we described in Lemma \ref{lem:FactorizationQSt}. This shows that to give a morphism between dagger Stein algebras is the same as to give a morphism of dagger Stein spaces and vice versa.
\ \hfill $\Box$

The following corollary is our version of the classical Forster's theorem. A similar statement holds for the category of affinoid algebras or 
in the dagger affinoid setting as discussed in Remark 4.18 of \cite{BaBe}.

\begin{cor} \label{cor:forster}
The functor defined by $\mathcal{M}$ from the opposite category of $k$--dagger Stein algebras to the category of locally ringed Grothendieck topological spaces is fully faithful. Therefore, we can identify the  category opposite to $k$-dagger Stein algebras, 
with a full sub-category of the category of locally ringed Grothendieck topological spaces. 
\end{cor}

Thanks to Theorem \ref{prop:quasi_Stein_algebra_spaces} we can give the following definition.

\begin{defn} \label{defn:stein_localization}
A morphism $A \to B$ between dagger Stein algebras is called a \emph{localization} if it can be written as a projective limit of dagger localizations.
\end{defn}

\begin{defn} \label{defn:stein_open_immersion}
	A morphism $f: X \to Y$ between dagger Stein spaces is called an \emph{open immersion} if it is a homeomorphism of $X$ with $f(X)$ and for each $x \in X$ the induced morphism of stalks $\mO_{Y, f(x)} \to \mO_{X, x}$ is an isomorphism.
\end{defn}

\begin{rem}
	Using Theorem \ref{prop:quasi_Stein_algebra_spaces} and Proposition 4.22 of \cite{BaBe}, one can see that $A \to B$ is a localization of Stein algebras if and only if the associated morphism of dagger Stein spaces $\mM(B) \to \mM(A)$ is an open immersion.
\end{rem}

\section{Homological characterization of the Stein topology} \label{sec:Stein_geometry}

We now need to prove last lemmata for deducing the main results of this work.

\begin{lem} \label{lemma:resolution}
Let $A$ be an object of $\ttComm(\ttCBorn_{k})$ and $\mathscr{L}_{A}^{\bullet}(E)$ be the Bar resolution (see Definition \ref{defn:Bar}) of an object $E\in \ttMod(A)$. If both $A$ and $E$ are nuclear, as object of $\ttCBorn_{k}$, then $\mathscr{L}_{A}^{\bullet}(E)$ is isomorphic to $E$ in $D^{\leq 0}(A)$.
\end{lem}
{\bf Proof.}
It is enough to check that we can apply Lemma \ref{lemma:flat_res}. This is possible because, by Theorem \ref{thm:strct_exact_nuclear}, nuclear objects of $\ttCBorn_{k}$ are flat for $\wotimes_k$.
\ \hfill $\Box$

Let $A \in \ttComm(\ttCBorn_k)$ and $E, F \in \ttMod(A)$. In Definition \ref{defn:Bar} we introduced the Bar resolution of $F$ as the complex $\cdots \to \mathscr{L}_{A}^{2}(F) \to \mathscr{L}_{A}^{1}(F) \to \mathscr{L}_{A}^0(F) \to F \to 0$
\[ \mathscr{L}_{A}^{n}(F) = A \wotimes_{k} (A^{\wotimes^{n}_{k}} \wotimes_{k} F) = A \wotimes_k (\underbrace{A \wotimes_k \cdots \wotimes_k A}_{n \text{ times}} \wotimes_k F) \]
whose differentials are given in equation (\ref{eqn:diff_bar}). Applying to this complex the functor $E \wotimes_A (-)$ one obtains the complex
\[ \cdots \to E \wotimes_A \mathscr{L}_{A}^{2}(F) \to E \wotimes_A \mathscr{L}_{A}^{1}(F) \to E \wotimes_A \mathscr{L}_{A}^{0}(F) \to E \wotimes_A F \to 0 \]
which is a representative of $E \wotimes_A^\Lb F$ in $D^{\le 0}(A)$.

\begin{rem} \label{rmk:bar_complex}
The terms of this complex can be written in a simplified form as follows
\[  E \wotimes_A \mathscr{L}_{A}^{n}(F) =  E \wotimes_A A \wotimes_k (A^{\wotimes^{n}_{k}} \wotimes_k F) \cong E \wotimes_k A^{\wotimes^{n}_{k}} \wotimes_k F. \]
\end{rem}

Therefore, we will introduce the notation
\begin{equation} \label{eqn:bar_notation}
\mathscr{L}_{A}^{n}(E, F) = E \wotimes_k A^{\wotimes^{n}_{k}} \wotimes_k F.
\end{equation}

\begin{lem} \label{lemma:new}
Let $A$, $B$ be dagger Stein algebras and let $X = \mM(B)$, $Y = \mM(A)$ be the corresponding dagger Stein spaces. Let $f: A \to B$ be a morphism and $\mM(f): X = \mM(B) \to Y = \mM(A)$ the corresponding morphism of dagger Stein spaces. Let also $C$ be another dagger Stein algebra and $g:A \to C$ an arbitrary morphism in $\ttComm(\ttCBorn_{k})$ Then, if $\mM(f)$ is an open immersion there exist projective systems of Banach algebras $\{ A_i \}_{i \in \Nb}$,  $\{ B_i \}_{i \in \Nb}$ and $\{ C_i \}_{i \in \Nb}$ and morphisms of systems $\{ f_i: A_i \to B_i \}_{i \in \Nb}$, $\{ g_i:A_i \to C_i \}_{i \in \Nb}$ such that
\begin{enumerate}
\item $\underset{i \in \Nb}\limpro A_i \cong A$, $\underset{i \in \Nb} \limpro B_i \cong B$, $\underset{i \in \Nb} \limpro C_i \cong C$ and the projection maps $A \to A_i$, $B \to B_i$, $C \to C_i$ have dense images;
\item $\limpro f_i = f$, $\limpro g_i = g$;
\item for every $i \in \Nb$ an isomorphism
\begin{equation} \label{eqn:der_iso}
 \mathscr{L}_{A_i}^{\bullet}(C_i, B_i) \to C_i \wotimes_{A_i} B_i 
 \end{equation} 
of objects of $D^{\le 0}(B_i)$.
\end{enumerate}
\end{lem}

{\bf Proof.}
By Lemma \ref{lem:FactorizationQSt} every morphism of Stein algebras can be represented as a system morphism of dagger affinoid algebras. So, let
\[ X = \bigcup_{i \in \Nb} U_i, \ \ Y = \bigcup_{i \in \Nb} V_i, \ \ Z = \mM(C) = \bigcup_{i \in \Nb} W_i \]
be representations of $X, Y$ and $Z$, with $U_i = \mM(B_i'), V_i = \mM(A_i')$ and $W_i= \mM(C_i')$ where $A'_i$, $B_i'$ and $C_i'$ are dagger affinoid algebras and the systems are chosen in a way that both morphisms $f$ and $A \to C$ have the same index set (for finite diagram without loops this re-indexing can always be performed; for more details see Proposition 3.3 of the Appendix of \cite{ArMaz}). Moreover, using Lemma \ref{lemma:Stein_algebra_strict} we can suppose that all $A'_i$, $B_i'$ and $C_i'$ are strictly dagger affinoid algebras.
Consider the morphisms $f_i': A_i' \to B_i'$ and $g_i: A_i' \to C_i'$, induced by $f$ and $g: A \to C$ respectively. It is easy to check that $\mM(f_i')$ are open immersions of dagger affinoid spaces.
Applying Proposition 4.22 of \cite{BaBe} we can deduce that $f_i'$ are dagger affinoid localizations. Applying Theorem 5.7 of \cite{BaBe}, we obtain a strict isomorphism of complexes of $B_i'$-modules 
\begin{equation} \label{eqn:iso_hom_epi}
C_i' \wotimes_{A_i'}^\Lb B_i' \cong C_i' \wotimes_{A_i'} B_i'.
\end{equation} 
The object $C_i' \wotimes_{A_i'}^\Lb B_i'$ can be represented by the complex $\mathscr{L}_{A_i'}^{\bullet}(C_i', B_i')$.
\begin{equation} \label{eqn:iso_hom_epi2} 
\cdots \to \mathscr{L}_{A_i'}^{2}(C_i', B_i') \to \mathscr{L}_{A_i'}^{1}(C_i', B_i') \to \mathscr{L}_{A_i'}^{0}(C_i', B_i') 
\end{equation} 
using notation of equation (\ref{eqn:bar_notation}). It is clear that the explicit description of the differentials of the Bar resolution given in equation (\ref{eqn:diff_bar}) implies that the Bar resolution of $B_i'$ admits a uniform indexing as a diagram of ind-Banach modules, once one fixes a representation of $B_i'$ as an ind-Banch object. Moreover, by Theorem 4.9 (4) of \cite{BaBe}  morphisms of dagger affinoid algebras can always be written as morphisms of inductive systems of Banach algebras.  Therefore, we can find representations of $A_i' \cong \underset{\rho > 1}\limind (A_\rho'')_i$ and $B_i' \cong \underset{\rho > 1}\limind (B_\rho'')_i$, $C_i' \cong \underset{\rho > 1}\limind (C_\rho'')_i$ such that the isomorphism (\ref{eqn:iso_hom_epi}) lifts to an isomorphism of complexes Banach modules for any $\rho > 1$ small enough. Indeed to choose such a $\rho$ we can proceed as follows. Let $\rho' > 1$ be small enough such that $f_i': A_i' \to B_i'$ is representable as a morphism of inductive systems $(A''_{\rho'})_i \to (B'_{\rho'})_i$ and let $\rho'' > 1$ be such that $g_i: A_i' \to C_i'$ can be represented as a morphism of inductive systems $(A'_{\rho''})_i \to (C'_{\rho''})_i$. Then, any $\rho$ such that $1 <\rho < \min \{ \rho', \rho''\}$ is a suitable choice. In fact, all the maps of the complex of equation (\ref{eqn:iso_hom_epi2}) can be written as a morphism of inductive systems for such small $\rho$, because these maps are obtained from the differentials of the Bar resolution by tensoring $(B_\rho'')_i$ (which becomes an $(A_\rho'')_i$-module) over with $(C_\rho'')_i$ (which becomes a $(A_\rho'')_i$-modules through the map $(A''_{\rho''})_i \to (C''_{\rho''})_i$). Therefore, we have a quasi-isomorphism of complexes of Banach modules
\begin{equation} \label{eqn:iso_hom_epi_ban}
[\cdots \to \mathscr{L}_{(A''_\rho)_i }^{2}((C''_\rho)_i , (B''_\rho)_i ) \to \mathscr{L}_{(A''_\rho)_i }^{1}((C''_\rho)_i , (B''_\rho)_i ) \to \mathscr{L}_{(A''_\rho)_i }^{0}((C''_\rho)_i , (B''_\rho)_i ) ]\to (C''_\rho)_i \wotimes_{(A''_\rho)_i } (B''_\rho)_i 
\end{equation}  
for $\rho$ small enough. Moreover, by the fact that the dagger subdomain embeddings
\[ U_i \subset U_{i + 1}, \ \ V_i \subset V_{i + 1}, \ \ W_i \subset W_{i + 1} \]
are inner for all $i \in \Nb$, we can find $\rho_U > 1$, $\rho_V > 1$, $\rho_W > 1$ small enough
 such that 
\[ \mM((B_{\rho_U}'')_i) \subset U_{i + 1}, \ \mM((A_{\rho_V}'')_i) \subset V_{i + 1}, \mM((C_{\rho_W}'')_i) \subset W_{i + 1}. \]
Therefore, fixing any $\rho > 1$ such that $ \min \{ \rho_U, \rho_V, \rho_W, \rho', \rho'' \} > \rho > 1$, we set for each $i$
\[ A_i = (A_\rho'')_i, \ \ B_i = (B_\rho'')_i, \ \ C_i = (C_\rho'')_i. \]
The properties (1), (2) and (3) satisfied by the system so obtained as direct consequences of the choices we made.
 \hfill $\Box$
 
 \begin{rem}
 Notice that the isomorphism (\ref{eqn:der_iso}) does not imply that $C_i \wotimes_{A_i}^\Lb B_i \to C_i \wotimes_{A_i} B_i$ is an isomorphism in general. Because for Archimedean base fields Banach spaces might not be flat in $\ttCBorn_k$ and therefore the Bar resolution is not a flat resolution.
 \end{rem}

\begin{thm} \label{thm_stein_homotopy}
Let $i: U \to V$ be an open immersion of dagger Stein spaces corresponding to a morphism $A_V\to B_U$ of dagger Stein algebras over $k$. For any dagger Stein $A_V$-algebra $C_W$ the morphism
\[ C_W \wotimes_{A_V}^\Lb B_U \to C_W \wotimes_{A_V} B_U \]
is an isomorphism in $D^{\le 0}(B_U)$. Also (in the case that $C_W = B_U$) $A_{V} \to B_{U}$ is a homotopy epimorphism.
\end{thm}
{\bf Proof.}
Lemma \ref{lemma:new}, applied to $A_V, B_U, C_W$ and the morphism $i: U \to V$ yields three systems of Banach algebras such that
\[ \underset{i \in \Nb} \limpro A_i \cong A_V, \ \ \underset{i \in \Nb} \limpro B_i \cong B_U, \ \ \underset{i \in \Nb} \limpro C_i \cong C_W \]
such that the opposite morphisms to $i$ and the morphism $A_V \to C_W$ can be written as morphisms of these projective systems. Moreover, for each $i \in \Nb$, Lemma \ref{lemma:new} also yields to a strictly exact complex
\[  \cdots \to \mathscr{L}_{A_i}^{2}(C_i, B_i) \to \mathscr{L}_{A_i}^{1}(C_i, B_i) \to \mathscr{L}_{A_i}^{0}(C_i, B_i) \to C_i \wotimes_{A_i} B_i \to 0 \]
The following commutative diagram
\[
\begin{tikzpicture}
\matrix(m)[matrix of math nodes,
row sep=2.6em, column sep=2.8em,
text height=1.5ex, text depth=0.25ex]
{  B_j \wotimes_k C_j & B_i \wotimes_k C_i  \\
   B_j \wotimes_{A_j} C_j & B_i \wotimes_{A_i} C_i   \\};
\path[->,font=\scriptsize]
(m-1-1) edge node[auto] {$$} (m-1-2);
\path[->,font=\scriptsize]
(m-1-1) edge node[auto] {$$} (m-2-1);
\path[->,font=\scriptsize]
(m-1-2) edge node[auto] {$$}  (m-2-2);
\path[->,font=\scriptsize]
(m-2-1) edge node[auto] {$$}  (m-2-2);
\end{tikzpicture}
\]
show that the bottom horizontal map is an epimorphism because the rest of the arrows are. Continuing to tensor we get a epimorphism $\mathscr{L}_{A_j}^{n}(C_j, B_j) \to \mathscr{L}_{A_i}^{n}(C_i, B_i)$ for any $n$. It is easy to check that these morphisms commutes with differentials, giving a projective system of morphism of complexes
\begin{equation} \label{eqn:system_complex}
 \begin{tikzpicture}
 \matrix(m)[matrix of math nodes,
 row sep=2.6em, column sep=2.8em,
 text height=1.5ex, text depth=0.25ex]
 {  
 	\cdots &\mathscr{L}_{A_j}^{1}(C_j, B_j) & \mathscr{L}_{A_j}^{0}(C_j, B_j) & C_j \wotimes_{A_j} B_j & 0  \\
    \cdots &\mathscr{L}_{A_i}^{1}(C_i, B_i) & \mathscr{L}_{A_i}^{0}(C_i, B_i) & C_i \wotimes_{A_i} B_i & 0  \\};
 \path[->,font=\scriptsize]
 (m-1-1) edge node[auto] {$d_j^2$} (m-1-2);
 \path[->,font=\scriptsize]
 (m-1-2) edge node[auto] {$d_j^1$} (m-1-3);
 \path[->,font=\scriptsize]
 (m-1-3) edge node[auto] {$d_j^0$} (m-1-4);
 \path[->,font=\scriptsize]
 (m-1-4) edge node[auto] {$$} (m-1-5);
 \path[->,font=\scriptsize]
 (m-1-2) edge node[auto] {$$} (m-2-2);
 \path[->,font=\scriptsize]
 (m-1-3) edge node[auto] {$$} (m-2-3);
 \path[->,font=\scriptsize]
 (m-1-4) edge node[auto] {$$} (m-2-4);
 \path[->,font=\scriptsize]
 (m-2-1) edge node[auto] {$d_i^2$} (m-2-2);
 \path[->,font=\scriptsize]
 (m-2-2) edge node[auto] {$d_i^1$} (m-2-3);
 \path[->,font=\scriptsize]
 (m-2-3) edge node[auto] {$d_i^0$} (m-2-4);
 \path[->,font=\scriptsize]
 (m-2-4) edge node[auto] {$$} (m-2-5);
 \end{tikzpicture}
 \end{equation}
of whose limit we denote by
\begin{equation} \label{eqn:limit_complex}
\cdots \longrightarrow \limpro_{i \in \Nb} (C_{i}\wotimes_{k} A_{i} \wotimes_{k} B_i) \longrightarrow \limpro_{i \in \Nb} (C_{i}\wotimes_{k} B_i) \longrightarrow \limpro_{i \in \Nb} (C_{i}\wotimes_{A_i} B_i) \longrightarrow 0.
\end{equation}
Using the isomorphisms  $\coim(d_i^{n+1}) \cong \ker (d_i^n)$ and $\coim(d_j^{n+1}) \cong \ker(d_j^n)$ the commutative squares
\begin{equation}
\xymatrix{
\mathscr{L}_{A_j}^{n}(C_j, B_j)\ar[r] \ar[d] & \ker(d_j^{n-1}) \ar[d] \\
\mathscr{L}_{A_i}^{n}(C_i, B_i)\ar[r] & \ker(d_i^{n-1}) 
}
\end{equation}
are seen to have all arrows epimorphisms except the one on between the kernels. Therefore, it also is an epimorphism. Now both the systems $\{ \ker(d_i^{n}) \}_{i \in \mathbb{N}}$ and $\{\mathscr{L}_{A_i}^{n}(C_i, B_i)\}_{i \in \mathbb{N}}$ are epimorphic. Notice that 
\[\underset{i \in \mathbb{N}}\limpro  \ker(d_i^{n}) \cong \ker(\underset{i \in \mathbb{N}}\limpro d_{i}^{n}).
\]
Proposition \ref{prop:nuclear_proj} implies that $\underset{i \in \mathbb{N}}\limpro \mathscr{L}_{A_i}^{n}(C_i, B_i)$ is nuclear. Therefore, $\underset{i \in \mathbb{N}}\limpro  \ker(d_i^{n})$ is nuclear as it can be identified with a closed subspace of $\mathscr{L}_{A_V}^{n}(C_W, B_U)$. Therefore, Corollary \ref{cor:BornML2} implies that $\{ \ker(d_i^{n}) \}_{i \in \mathbb{N}}$ and $\{\mathscr{L}_{A_i}^{n}(C_i, B_i)\}_{i \in \mathbb{N}}$ are $\underset{i \in \mathbb{N}}\limpro$ acyclic.

Applying Lemma \ref{lem:LongShort} the to the system of complexes (\ref{eqn:system_complex}) we deduce that complex of equation (\ref{eqn:limit_complex}) is strictly acyclic. 
Moreover, applying Corollary \ref{cor:proj_lim_1} to each term in degree strictly less than zero of the complex (\ref{eqn:system_complex}), we see that (\ref{eqn:limit_complex}) is strictly isomorphic to the complex
\begin{equation} \label{eqn:limit_complex_2}
\cdots \longrightarrow (\limpro_{i \in \Nb} C_{i} )\wotimes_{k} (\limpro_{i \in \Nb} A_{i}) \wotimes_{k} (\limpro_{i \in \Nb} B_i) \longrightarrow ( \limpro_{i \in \Nb} C_{i} )\wotimes_{k} (\limpro_{i \in \Nb} B_i ) \longrightarrow \limpro_{i \in \Nb} (C_{i}\wotimes_{A_i} B_i) \longrightarrow 0.
\end{equation}
Thus we showed that the complex 
\[\cdots \longrightarrow  C_{W}\wotimes_{k} A_{V} \wotimes_{k} B_U \longrightarrow C_{W}\wotimes_{k} B_U \longrightarrow \limpro_{i \in \Nb} (C_{i}\wotimes_{A_i} B_i) \longrightarrow 0
\]
is strictly exact. Since $\coker(\mathscr{L}_{A_V}^{1}(C_W, B_U) \to \mathscr{L}_{A_V}^{0}(C_W, B_U)) \cong C_W \wotimes_{A_V} B_U$ we have shown the isomorphism in $D^{\le 0}(B_U)$
\[C_W \wotimes_{A_V}^\Lb B_U \to \limpro_{i \in \Nb} (C_{i}\wotimes_{A_i} B_i) \cong C_W \wotimes_{A_V} B_U.
\]
In the special case that the system $\{B_{i}\}_{i \in \Nb}$ and the system $\{C_{i}\}_{i \in \Nb}$ are equal we have $B_{U}=C_{W}$ and that $B_{i} \wotimes^{\mathbb{L}}_{A_i} B_{i} \cong B_{i} \wotimes_{A_i} B_{i} \to B_{i}$ is an isomorphism by Theorem 5.11 of \cite{BaBe}. We can infer that 
\[ B_U \wotimes_{A_V}^\Lb B_U \to  \limpro_{i \in \Nb}(B_{i} \wotimes_{A_{i}}B_{i}) \cong \limpro_{i \in \Nb}B_{i}  \cong B_{U} \]
is an isomorphism, concluding the proof of the theorem.
\ \hfill $\Box$

\begin{lem}\label{lemma:proj_homotopy_epi}
	Let $A$ be a dagger Stein algebra presented by a system of (strictly) dagger affinoid algebras $A_{i}$. Then the projection maps $A \to A_i$ are homotopy epimorphisms.
\end{lem}
{\bf Proof.}
We can write $A_i$ as a direct limit
\[ A_i \cong \limind_{\rho > 1} (A_i)_\rho \]
where $(A_i)_\rho$ are Stein algebras of Stein spaces that admit closed embeddings in polydisks and where each morphism $(A_i)_\rho \to (A_i)_{\rho'}$, for $\rho' < \rho$ corresponds geometrically to an open embedding. Such a system of Stein spaces can be found via a presentation of $A_i \cong  \frac{W_k^n}{I}$ and writing $W_k^n$ as the direct limit of the Stein algebras of open polydisks of radius bigger that one, which form a base of neighborhoods of the closed unital polydisk. Applying Theorem \ref{thm_stein_homotopy} we can infer that each $(A_i)_\rho \to (A_i)_{\rho'}$ is a homotopy epimorphism and applying Lemma \ref{lem:ind_limit_homotopy_epi} we can deduce that the canonical morphisms $(A_i)_\rho \to A_i$ are homotopy epimorphisms. Since $A$ is a Stein algebra and $A \to A_i$ corresponds geometrically to an open embedding there exists a $\rho > 1$ small enough such that $A \to A_i$ factors through $A \to A_\rho$ and that  $A \to A_\rho$ corresponds geometrically to an open embedding. Applying Theorem \ref{thm_stein_homotopy} we can deduce that $A \to A_\rho$ is a homotopy epimorphism and therefore $A \to A_i$ is a homotopy epimorphism because it can be written as a composition of two homotopy epimorphisms.
\ \hfill $\Box$

And finally, the last result of characterization of open Stein immersions.

\begin{thm}\label{thm:DaggerQStHoEpToImm}
Let $f: A_V \to B_U$ be a morphism of dagger Stein algebras over $k$. If $f$ is a homotopy epimorphism as morphism of $\ttComm(\ttCBorn_k)$, then it is a localization.
\end{thm}
{\bf Proof.}
The condition of $f$ being a homotopy epimorphism means that
\[ B_U \wotimes_{A_V}^\Lb B_U \cong B_U \wotimes_{A_V} B_U \cong B_{U}. \] 
Let $\underset{i \in \Nb} \limpro A_i \cong A_V, \underset{i \in \Nb} \limpro B_i \cong B_U$ be representations of $A_V$ and $B_U$ such that $f$ can be written as a morphism of projective systems $f_i:A_i \to B_i$ of dagger affinoid algebras. 
By Lemma \ref{lemma:proj_homotopy_epi} the projections $A_V \to A_i$ and $B_U \to B_i$ are homotopy epimorphisms. Lemma \ref{lem:composition_HomotopyMon} applied to the bottom horizontal map of the commutative diagram
\[
\begin{tikzpicture}
\matrix(m)[matrix of math nodes,
row sep=2.6em, column sep=2.8em,
text height=1.5ex, text depth=0.25ex]
{  A  & B  \\
   A_i & B_i   \\};
\path[->,font=\scriptsize]
(m-1-1) edge node[auto] {$$} (m-1-2);
\path[->,font=\scriptsize]
(m-1-1) edge node[auto] {$$} (m-2-1);
\path[->,font=\scriptsize]
(m-1-2) edge node[auto] {$$}  (m-2-2);
\path[->,font=\scriptsize]
(m-2-1) edge node[auto] {$$}  (m-2-2);
\end{tikzpicture}
\]
implies that $A_i \to B_i$ is a homotopy epimorphism for every $i \in \Nb$. This is equivalent to say that $f$ can be written as a projective system of homotopy epimorphisms of dagger affinoid algebras. Applying Theorem 5.11 of \cite{BaBe} we obtain that the morphisms $A_i \to B_i$ are open immersions of dagger affinoid spaces. Therefore, $\mM(f): \mM(B_U) \to \mM(A_V)$ can be written as filtered inductive limit of open embeddings, and it is easy to check that this implies that $\mM(f)$ is an open immersion.
\ \hfill $\Box$

To conclude our homological and categorical characterization of the topology of Stein spaces it remains to characterize coverings. Consider a Stein space $X$ and an arbitrary covering
\[ X = \bigcup_{i \in I} Y_i \]
of $X$ made of Stein spaces $Y_i$. By definition the topology of $X$ is hemi-compact and $X$ is also paracompact. Therefore, the family $\{ Y_i \}_{i \in I}$ admits a countable sub-family $\{ Y_j \}_{j \in J}$, where $J \subset I$, such that 
\[ X = \bigcup_{j \in J} Y_j. \]

\begin{lem} \label{lem:cover_one_way}
Let $A$ be a dagger Stein algebra and let $\{f_i: A \to A_{V_i} \}_{i\in I}$ be a family of localizations such that for some countable subset $J \subset I$ the corresponding family of functors
\[ \{F_{i}: \ttMod^{RR}(A) \to \ttMod^{RR}(A_{V_i})\}_{i \in J} \]
for is conservative. Then, the morphism $\phi: \underset{i \in J} \coprod \Max(A_{V_i} ) \to \Max(A)$ is surjective.
\end{lem}
{\bf Proof.}
Assume that the family of functors $\{F_i\}_{i \in J}$ is conservative and that $\phi: \underset{i \in J}\coprod \Max(A_{V_i} ) \to \Max(A)$ is not surjective. We will deduce a contradiction. Let $\fm_x \in \Max(A)$ be a point which is not in the image of $\phi$. Consider the quotient $A/\fm_x$. This is a Stein algebra: see Example \ref{exa:stein}. By \ref{thm_stein_homotopy} $A/\fm_x\in \ttMod^{RR}(A)$ and it is a non-trivial module. We also have  
\[ A_{V_i} \wotimes_A (A/\fm_x) = 0 \]
for all $i$, because the extension of $\fm_x$ to $A_{V_i}$ is equal to the improper ideal for all $i \in J$. This proves that the family $F_i$ is not conservative.
\ \hfill $\Box$

\begin{defn}\label{defn:RRqcoh_F}
We denote with $\ttMod_F^{RR}(A)$ the full sub-category of $\ttMod^{RR}(A)$ for which $M \in \ttMod^{RR}(A)$ is a bornological Fr\'echet space.
\end{defn}

\begin{lem} \label{prop:proj_lim_A}
	Let $\{ E_i \}_{i\in I}$ be a countable set of binuclear bornological Fr\'echet modules over $A$ and $F$ a bornological Fr\'echet modules over $A$, where $A$ is a nuclear Fr\'echet algebra. Then, the canonical map
	\[ E \wotimes_{A} F = (\prod_{i\in I} E_i) 
	\wotimes_{A} F \to \prod_{i \in I}(E_i \wotimes_{A} F)  \] 
	is an isomorphism of bornological modules.
\end{lem}
{\bf Proof.}
\[ \prod_{i \in I} F \wotimes_A E_i = \prod_{i \in I} \coker ( F \wotimes_k A \wotimes_k E_i \stackrel{d_1}\to F \wotimes_k E_i) \]
where $d_1$ is induced by the differential in degree $1$ of the Bar resolution. Since direct products of bornological spaces preserve cokernels (see Proposition 1.9 of \cite{PrSc} for a proof of this fact for bornological spaces over $\Cb$ and Proposition 1.2.12 of \cite{Bam} for the same proof worked out in a more general settings) we see that
\[ \prod_{i \in I} \coker ( F \wotimes_k A \wotimes_k E_i \stackrel{d_1}\to F \wotimes_k E_i) \cong  \coker ( \prod_{i \in I} ( F \wotimes_k A \wotimes_k E_i) \stackrel{d_1}\to \prod_{i \in I} (F \wotimes_k E_i)) \] 
to which we can apply Corollary \ref{cor:proj_lim_2} (cofiltering the infinite direct product by its finite products) to deduce that
\[ \coker ( \prod_{i \in I} (F \wotimes_k A \wotimes_k E_i) \stackrel{d_1}\to \prod_{i \in I} (F \wotimes_k E_i)) \cong \coker ( F \wotimes_k A \wotimes_k (\prod_{i \in I} E_i) \stackrel{d_1}\to F \wotimes_k (\prod_{i \in I} E_i)) \cong F \wotimes_{A} E.  \] 
\ \hfill $\Box$

\begin{cor} \label{cor:proj_lim_A}
	Under the same hypothesis of Lemma \ref{prop:proj_lim_A} we have that for any countable set $I$,
	\[ E \wotimes_{A}^\Lb F = (\prod_{i\in I} E_i) 
	\wotimes_{A}^\Lb F \to \prod_{i \in I}(E_i \wotimes_{A}^\Lb F)  \] 
	is an isomorphism of bornological modules.
\end{cor}
{\bf Proof.}
The same reasoning of Lemma \ref{prop:proj_lim_A} can be extended to $E \wotimes_{A}^\Lb F$ representing it with the complex $\sL_A^\bullet(E, F)$, using the notation introduced so far. Therefore, for each $n \in \Nb$ we have that
\[ \sL_A^n(E, F) = E \wotimes_k A^{\wotimes n} \wotimes_k F = (\prod_{i \in I} E_i )\wotimes_k A^{\wotimes n} \wotimes_k F. \] When $I$ is finite, using the fact that the completed projective tensor product commutes with finite products, we can deduce that $\sL_A^n(E, F) \cong \underset{i \in I}\prod (E_i \wotimes_k A^{\wotimes n} \wotimes_k F)$. By writing a coundable product as a cofiltered projective system of finite products, we can use Corollary \ref{cor:proj_lim_2} to deduce that
\[ \sL_A^n(E, F) \cong \prod_{i \in I} (E_i \wotimes_k A^{\wotimes n} \wotimes_k F) \]
for general countable collections. Now $\underset{i \in I} \prod$ defines a functor $\underset{i \in I} \prod D^{\le 0}(\ttMod(A)) \to D^{\le 0}(\ttMod(A))$ because direct products are exact and so we get a strict quasi-isomorphism of the complexes representing $E \wotimes_{A}^\Lb F $ and $\underset{i\in I}\prod(E_i \wotimes_{A}^\Lb F)$.

\ \hfill $\Box$

\begin{lem} \label{lem:cover_the_other_way}
	Let $A$ be a Stein algebra and $\{ V_i \}_{i \in \Nb}$ a countable collection of Stein domains that covers $X = \Max(A)$. Then, the corresponding family of functors $\ttMod_F^{RR}(A) \to \ttMod_F^{RR}(A_{V_i})$ is conservative.
\end{lem}
{\bf Proof.}
Let $f: M \to N$ in $\ttMod_F^{RR}(A)$ be any morphism such that $f_{i}: M \wotimes_{A}A_{V_i} \to N \wotimes_{A}A_{V_i} $ are isomorphisms for all $i$. The \v{C}ech-Amitsur complex
\begin{equation} \label{eqn:cech_amistur}
 0 \to A \to \prod_{i_1 \in \Nb} A_{V_{i_1}} \to \prod_{i_1, i_2 \in \Nb} A_{V_{i_1}} \wotimes_{A} A_{V_{i_2}} \to \cdots 
\end{equation}
is strictly exact as a consequence of the Theorem B for Stein spaces (cf. Fundamental Theorem in page 124 of \cite{GR} for the Archimedean version of Theorem B and Satz 2.4 \cite{KI} for the non-Archimedean one). Theorem \ref{thm:extend_tensor} permits to apply the functor $M \wotimes_{A}^{\mathbb{L}}(-)$ to the complex (\ref{eqn:cech_amistur}) because it permits to extend $M \wotimes_{A}^{\mathbb{L}}(-)$ to a functor on the derived category of unbounded complexes. Therefore, by applying the derived functor $M \wotimes_{A}^{\mathbb{L}}(-)$ we are left with an object strictly quasi-isomorphic to zero. Notice that in this case that $M \wotimes^{\mathbb{L}}_{A} (-)$ commutes with the relavant countable products (by Corollary \ref{cor:proj_lim_A}) because $M$ is Fr\'echet and products respect cokernels. Furthermore,  as $M$ is RR-quasi-coherent, we have \[M \wotimes^{\mathbb{L}}_{A} (A_{V_{i_1}}\wotimes_{A} \cdots \wotimes_{A} A_{V_{i_n}} ) \cong M \wotimes_{A} (A_{V_{i_1}}\wotimes_{A} \cdots \wotimes_{A} A_{V_{i_n}} )  \]
for each $i_i, \dots i_n$. Therefore, we can apply Corollary \ref{cor:proj_lim_2} and by a small computation, obtain a a strictly exact complex
\[0 \to M \to  \prod_{i_1 \in \Nb} (M \wotimes_{A} A_{V_{i_1}}) \to \prod_{i_1,i_2 \in \Nb} (M \wotimes_{A} A_{V_{i_1}} \wotimes_{A} A_{V_{i_2}}) \to \cdots
\]
We can do the same thing for $N$ and this yields that the morphisms $f_{i}$ extend uniquely to morphisms of the complexes resolving $M$ and $N$. Therefore $f$ is an isomorphism.
\ \hfill $\Box$

\begin{thm} \label{thm:coverings}
A collection of homotopy epimorphisms of Stein algebras $\{A_V \to A_{V_i}\}_{i \in I}$ is a covering of $\Max(A_V)$ if and only if the family of functors $\{\ttMod_F^{RR}(A) \to \ttMod_F^{RR}(A_{V_i})\}_{i \in J}$ is conservative, with $J \subset I$ a countable subset.
\end{thm}
{\bf Proof.}
The assertion of the theorem is simply the combination of Lemma \ref{lem:cover_one_way} and Lemma \ref{lem:cover_the_other_way}.
\ \hfill $\Box$

We summarize the main results of this section in the following corollary.

\begin{cor} \label{cor:main_results}
	Consider $(\ttCBorn_k, \wotimes_k, k)$ as a closed symmetric monoidal elementary quasi-abelian category. The natural inclusion of categories $\ttStn_k \rhook \ttComm(\ttCBorn_k)$ permits to define a countable version of the formal homotopy Zariski topology on $\ttStn_k$ as in Definition \ref{defn:homotopy_Zariski}. The coverings of Stein spaces by Stein spaces (in the usual sense for Archimedean base fields, and in the rigid sense for non-Archimedean base fields) corresponds precisely with the families of morphisms $\{A_V \to A_{V_i}\}_{i \in I}$ in the category $\ttComm(\ttCBorn_{k})$ for which the family of functors $\{\ttMod_F^{RR}(A) \to \ttMod_F^{RR}(A_{V_i})\}_{i \in J}$ is conservative, and $A_V \to A_{V_i}$ is a homotopy epimorphism for each $i \in J$ with $J \subset I$ some countable subset.
\end{cor}
{\bf Proof.}
Theorems \ref{thm:DaggerQStHoEpToImm} and \ref{thm_stein_homotopy} precisely mean that in $\ttStn_k$ a morphism is a homotopy epimorphism if and only if it is an open immersion. Whereas, the claim on the coverings is obtained in Theorem \ref{thm:coverings}.
\ \hfill $\Box$

\begin{rem}If instead of allowing infinite covers, we want to consider only finite covers of Stein spaces by Stein spaces (so that $J$ is finite), then the analogue of Corollary \ref{cor:main_results} gives a description of the formal homotopy Zariski topology and in this case one can replace $\ttMod_F^{RR}$ with $\ttMod^{RR}$, using all RR-quasicoherent modules instead of just the Frechet ones. Notice also that within the proof of Lemma \ref{lem:cover_the_other_way} we have shown that given a countable open Stein cover of a Stein space and any $M \in \ttMod_F^{RR}(A)$ that the complex 
\[0 \to M \to  \prod_{i_1 \in \Nb} (M \wotimes_{A} A_{V_{i_1}}) \to \prod_{i_1,i_2 \in \Nb} (M \wotimes_{A} A_{V_{i_1}} \wotimes_{A} A_{V_{i_2}}) \to \cdots
\]
is strictly exact. The same holds for finite open Stein covers of a Stein space and any $M \in \ttMod^{RR}(A)$.
\end{rem}

\section{Conclusions} 

It is natural to ask if the characterization of the immersions of dagger Stein spaces given by Theorems \ref{thm_stein_homotopy} and \ref{thm:DaggerQStHoEpToImm} can be extended to more general kind of domains. As a first attempt to generalizations one can consider quasi-Stein algebras and quasi-Stein spaces in the sense of \cite{KI}, Definition 2.3. Therefore, in the setting of classical (\ie non-dagger) non-Archimedean geometry, one can define quasi-Stein spaces as analytic spaces which have an exhaustion analogous to the one we described in Definitions \ref{defn:stein_algebra} and \ref{defn:stein_algebra}, with the difference that the immersion $U_i \to U_{i + 1}$ can be any Weierstrass embedding. More precisely, for a non-Archimedean non-trivially valued base field the following definition makes sense.

\begin{defn} \label{defn:quasi_stein_algebra}
	A \emph{quasi-Stein algebra} over $k$ is a complete bornological algebra $A$ over $k$ which is isomorphic to an inverse limit of $k$ affinoid algebras 
\begin{equation} \label{eqn:quasi_stein}
	\cdots \longrightarrow A_4 \longrightarrow A_3 \longrightarrow A_2 \longrightarrow A_1 \longrightarrow A_0
	\end{equation} 
	in the category $\ttCBorn_{k}$ where each morphism is a Weierstrass localization for each $i$.
\end{defn}

We notice that the quasi-Stein algebras that one obtains in this way turn out to be Fr\'echet algebras. Therefore, the theory developed so far is meaningful also for these kind of algebras and their associated quasi-Stein spaces, and all the arguments work proving also in that case that open immersions can be characterized homologically as we did for Stein spaces. There is only one small difference to take into account. Since the Weierstrass localization of equation (\ref{eqn:quasi_stein}) are not required to be inner, they are not nuclear map of the underlying Banach spaces and for this reason, quasi-Stein algebras fail to be nuclear. This implies that the bornological version of Mittag-Leffler Lemma, see Lemma \ref{lemma:mittag_frechet}, does not apply. However, one can deduce the results of section \ref{sec:Stein_geometry} for quasi-Stein algebras by working with the compactoid bornology instead of the von Newmann bornology. We do not know if the same results hold when quasi-Stein algebras are equipped with the von Neumann bornology.

The dagger version of the notion of quasi-Stein algebras and quasi-Stein spaces can be defined by analogy and that same questions can be asked for them. This case is much more complicated and it is even possible that the analogous of Theorems \ref{thm_stein_homotopy} and \ref{thm:DaggerQStHoEpToImm} for dagger quasi-Stein spaces does not hold.

As concluding remarks, we spend few words to describe an outlook  for potential future work. In \cite{BeKr2} a basic framework for derived analytic geometry will be given and the results in this paper will then explain the relation between the derived analytic spaces and usual analytic spaces. Also, the notion of RR-module, introduced in Definition \ref{defn:RRqcoh}, can be further developed and used to define quasi-coherent sheaves on analytic spaces. This notion will be useful to give a more general definition of derived analytic spaces with compared with the ones proposed in \cite{Por} (Definition 1.3) and \cite{PY} (Defintion 2.5), where one obstacle to mimic the definitions of derived scheme is the lack of a good notion of quasi-coherent sheaf in analytic geometry.

\bibliographystyle{amsalpha}

\end{document}